\documentclass[a4paper, 11pt]{article}
\pdfoutput=1
\usepackage[a4paper, top=30mm, bottom=30mm, left=30mm, right=30mm]{geometry}
\usepackage{setspace} 
\usepackage[T1]{fontenc}
\usepackage[utf8]{inputenc}
\usepackage{lmodern}        
\usepackage[english]{babel} 
\usepackage[hyphens]{url}
\usepackage{hyperref}
\usepackage{amsmath}
\usepackage{mathtools}
\usepackage{amssymb}
\usepackage{amsfonts}
\usepackage{bm}
\usepackage{amsthm}
\usepackage[numbers]{natbib}
\usepackage{float}
\usepackage{graphicx}
\usepackage{booktabs}
\usepackage{adjustbox}
\usepackage{rotating}
\usepackage{subfig}
\usepackage[font=small, labelfont=bf]{caption}
\usepackage{datetime}
\usepackage{calc}
\usepackage{xcolor, color}
\usepackage{xspace}
\usepackage{enumitem}
\usepackage{soul}
\usepackage{emptypage} 
\usepackage{cleveref}

\hypersetup{%
pdfauthor={},%
pdftitle={},%
pdfsubject={Subject},%
pdfkeywords={Keywords},%
pdfproducer={}, pdfcreator={pdfTeX - Emacs - Linux },%
colorlinks=true, breaklinks=true, %
citecolor=blue, linkcolor=blue, urlcolor=magenta, %
bookmarksnumbered=true,bookmarksopen=true,bookmarksopenlevel=1,%
pdftex=true%
}%

\DeclareGraphicsExtensions{.pdf,.png,.jpg,.jpeg}

\definecolor{colNOTE}{HTML}{ED2E38}

\newcommand*{\eg}{e.g.\xspace}
\newcommand*{\ie}{i.e.\xspace}

\makeatletter
\newcommand*{\etc}{%
\@ifnextchar{.}%
{etc}%
{etc.\@\xspace}%
}
\makeatother

\newcommand{\figref}[1]{Figure~\ref{#1}}

\newcommand{\nst}{{\raise0.5ex\hbox{\footnotesize{st}}}\xspace}
\newcommand{\nnd}{{\raise0.5ex\hbox{\footnotesize{nd}}}\xspace}
\newcommand{\nrd}{{\raise0.5ex\hbox{\footnotesize{rd}}}\xspace}
\newcommand{\nth}{{\raise0.5ex\hbox{\footnotesize{th}}}\xspace}

\theoremstyle{plain}
\newtheorem{theorem}{Theorem}
\newtheorem{corollary}{Corollary}
\newtheorem{lemma}{Lemma}
\newtheorem{proposition}{Proposition}
\theoremstyle{definition}

\newtheorem{assumption}{Assumption}
\theoremstyle{remark}
\newtheorem{remark}{Remark}

\crefname{theorem}{theorem}{theorems}
\crefname{corollary}{corollary}{corollaries}
\crefname{lemma}{lemma}{lemmas}
\crefname{proposition}{proposition}{propositions}
\crefname{definition}{definition}{definitions}
\crefname{assumption}{assumption}{assumptions}
\crefname{remark}{remark}{remarks}
\crefname{example}{example}{examples}

\newlist{cnd}{enumerate}{1} 
\setlist[cnd,1]{label={\upshape (C\arabic*)}, ref={\upshape (C\arabic*)},
                topsep=0ex,itemsep=0.5ex,partopsep=0pt, parsep=0pt,
                , resume} 
\crefname{cndi}{condition}{conditions} 

\newcommand{\eps}{\varepsilon}

\newcommand{\ind}{\mathbb{I}}
\newcommand{\R}{\mathbb{R}}
\newcommand{\Rpos}{\R_{>0}} 
\newcommand{\Rnneg}{\R_{\geq 0}} 
\newcommand{\N}{\mathbb{N}}

\newcommand{\set}[1]{\left\{ #1 \right\}}
\newcommand{\abs}[1]{\left\lvert #1 \right\rvert}
\newcommand{\norm}[1]{\left\lVert #1 \right\rVert}

\newcommand{\inv}[1]{#1\raisebox{1.15ex}{$\scriptscriptstyle-\!1$}}
\newcommand{\tr}[1]{#1\raisebox{1.15ex}{$\scriptscriptstyle\!\textsf{T}$}}

\DeclareMathOperator*{\argmax}{arg\,max}

\DeclareMathOperator{\oPr}{Pr}
\newcommand{\prob}[2][]{{\oPr_{#1}}{\left[#2\right]}}
\DeclareMathOperator{\oE}{E}
\newcommand{\E}[2][]{{\oE_{#1}}{\left[#2\right]}}
\DeclareMathOperator{\oVar}{Var}
\newcommand{\Var}[2][]{{\oVar_{#1}}{\left[#2\right]}}

\newcommand{\given}{\,\vert\,}

\newcommand{\conv}{\longrightarrow}

\newcommand{\bI}{\bm{I}}

\newcommand{\bSigma}{\bm{\Sigma}}
\newcommand{\bOmega}{\bm{\Omega}}

\newcommand{\bzero}{\bm{0}}
\newcommand{\bx}{\bm{x}}
\newcommand{\by}{\bm{y}}

\newcommand{\bmu}{\bm{\mu}}
\newcommand{\btheta}{\bm{\theta}}
\newcommand{\brho}{\bm{\rho}}

\newcommand{\bk}{\bm{k}}
\newcommand{\bkappa}{\bm{\kappa}}
\newcommand{\bdelta}{\bm{\delta}}

\newcommand{\lmin}{\lambda_{\min}}
\newcommand{\lmax}{\lambda_{\max}}

\DeclareMathOperator{\cl}{cl}

\DeclareMathOperator{\ESD}{ESD}

\newcommand{\REV}[1]{#1}

\begin{document}

\title{\MakeUppercase{\bf \normalsize %
Consistency for constrained maximum likelihood estimation and clustering based on mixtures of elliptically-symmetric distributions under general data generating processes
}}

\author{%
Pietro Coretto\thanks{%
Departiment of Economics and Statistics; %
University of Salerno (Italy) %
-- E-mail: \texttt{pcoretto@unisa.it}.
For this research the author was supported by Ministero dell'Università e della Ricerca (MUR, Italy), grant number   20223725WE.} %
\and
Christian Hennig\thanks{%
Department of Statistics; %
University of Bologna (Italy) %
--  E-mail: \texttt{christian.hennig@unibo.it}%
} 
}


\thispagestyle{empty}
\maketitle

\begin{abstract} \noindent
The consistency of the maximum likelihood estimator for mixtures of elliptically-symmetric distributions for estimating its population version is shown, where the underlying distribution $P$ is nonparametric and does not necessarily belong to the class of mixtures on which the estimator is based. In a situation where $P$ is a mixture of well enough separated but nonparametric distributions it is shown that the components of the population version of the estimator correspond to the well separated components of $P$. This provides some theoretical justification for the use of such estimators for cluster analysis in case that $P$ has well separated subpopulations even if these subpopulations differ from what the mixture model assumes.
\vspace{2ex}

\end{abstract}

\section{Introduction}\label{sec:intro}

Clustering methods are widely used in statistics, computer science, and many applied fields of science to discover group structures and to deal with the intrinsic inhomogeneity of complex data sets.
While there is abundant work on methodology, algorithms, and applications, a smaller body of literature has investigated the relationship between the clusters found by a method and the underlying data-generating mechanism.
Assuming that the observed data set is generated by independent and identical observations from a probability law $P$, consistency concerns the relationship between $P$ and the outcome of a method for random samples of a size converging to infinity. In cluster analysis, the clustering itself and/or distributional parameters characterising the clustering may be of interest.

Here we will derive consistency results for model-based clustering, i.e., clustering based on probability mixture models. More precisely, the results will concern maximum likelihood (ML) estimators (MLE) of finite mixtures of distributions from elliptically symmetrical distribution (ESD) families such as the Gaussian distribution. Finite mixture models (FMM) are convex combinations of probability distributions suitable to represent inhomogeneous populations. FMMs are a versatile tool for modelling complex distributions, and are at the basis of a variety of data analysis, prediction, and inference tasks: density estimation, regression, clustering, and classification \cite[see][]{fcr19-handbook,hmmr16-handbook}.

Given $K$, the number of mixture components, FMMs are represented by parameters characterising the individual mixture components (such as Gaussian mean and variance) and the component proportions. MLEs can be obtained in practice using the Expectation-Maximisation (EM) algorithm of \cite{delaru77}. Given the fitted parameters, a partition of the data can be obtained by use of the maximum posterior assignment rule, see \Cref{sec:mix-clustering-classification}.

In clustering, clusters are usually interpreted as corresponding to the estimated mixture components, although this is not always appropriate if different mixture components are not well separated, see \cite{hennig10}. 

FMMs are parametric, and therefore, as a standard, statisticians are interested in whether the parameters can be consistently estimated if the true underlying $P$ is indeed such a mixture. Consistency theory for ML estimation of mixture parameters is not as easy to obtain as for standard ML estimation, and standard conditions such as required by \cite{wald49} are not fulfilled. \cite{kw56} noted that the likelihood function of univariate Gaussian mixtures is unbounded if one of the mixture components' variances is allowed to converge to zero from above. This issue occurs for many classes of FMMs including those of ESDs as treated here.
There have been several proposals to solve it. 
A popular strategy, originally due to \cite{dennis81}, is to impose an ``eigen-ratio constraint'' (ERC). This is a constant that bounds the ratio between the largest and smallest eigenvalue across all components' scatter matrices.
The ERC allows for a non-compact parameter space; it is easier but also more restrictive to constrain the parameter space to be compact as done, e.g., by
\cite{redner81}, who proved consistency for the MLE of a general class of FMMs.
\cite{hathaway85} expanded the work of \cite{redner81} including the ERC in the case of a univariate Gaussian mixture. 
There is a long history of  ML-theory for mixtures. For reviews see \cite{RedWal84,chen17}.

Here, however, we focus on \REV{consistency in the case that} $P$ is general and not necessarily of the assumed mixture type. \REV{We will refer to this as ``nonparametric consistency'', but note that ``nonparametric clustering'' in the literature can refer to clustering methods that are not based on parametric model assumptions \cite{hunter24}, or to methods that do not treat the number of clusters as fixed \cite{rao16}. Here we investigate parametric clustering methods under nonparametric conditions.}

 Note that we will treat the number of estimated mixture components $K$ as fixed, so that the fitted mixture cannot approximate almost any $P$ by selecting a sufficiently large $K$. The results will apply to a $P$ that is indeed a mixture of the assumed family with $K$ mixture components, however the underlying $P$ could also have fewer or more than $K$ mixture components of this kind (in clustering it may sometimes be desirable to fit a mixture with more than $K$ not well separated components by a $K$-component mixture where the $K$ components are better separated and justified to be interpreted as clusters, see, e.g., \cite{BiCeGo00}), or it may not be possible to represent it as a mixture of the assumed type at all. 

In practice, arguably, parametric statistical model assumptions are never precisely fulfilled, yet often it is claimed that the application of data analytic methods based on parametric models require the parametric model assumptions to be fulfilled. In fact often standard theory relies on these model assumptions, and it may be suspected that if assumptions are violated, the method may not achieve the performance that is theoretically guaranteed assuming the model. Nonparametric theory as derived here, even if only asymptotic, will apply to the use of such methods more generally.

Nonparametric consistency results for cluster analysis go back to the seminal works of \cite{pollard81} (for $K$-means clustering) and \cite{hartigan81} (for single linkage). Further such results have been provided by \cite{ss12} for the DBSCAN algorithm and \cite{vbb08} for spectral clustering, among others. Multivariate Gaussian FMMs have been studied by  \cite{ggmm15} and \cite{ch17-jmlr}. The latter paper also has a nonparametric consistency result regarding a robust ML-type method accommodating outliers by incorporating an improper uniform mixture component covering the whole data space, and an Expectation-Conditional Maximization (ECM) algorithm implementing the ERC. 

There are two aspects of nonparametric consistency relevant in clustering (cp. \cite{vbb08}).
\begin{enumerate}[label=(\textit{\roman*})]
\item
Is the clustering method consistent at all, i.e., is there a limit clustering (or limit parameters characterising it) if the sample size grows larger and larger? Many clustering methods including MLE of FMMs are defined by parameters optimising an objective function. The asymptotic limit is then usually a functional defined on $P$ by generalizing the objective function to general distributions (``population version'' or ``canonical functional'') such as in \cite{pollard81,ggmm15,ch17-jmlr} and also in the present work. 

\item
Provided that a limit clustering exists, does this correspond to a reasonable clustering structure given $P$ that is of interest when doing cluster analysis?
\end{enumerate}

Addressing these aspects, the present paper provides two main contributions:
\begin{enumerate}[label=(\textit{\roman*})]

\item
We establish the existence and the consistency of the MLE for a general class of FMMs based on ESDs for a general nonparametric  $P$, generalizing the results for Gaussian FMMs in \cite{ggmm15,ch17-jmlr}. This implies finite sample existence of the MLE under mild conditions. 

\item
Under the additional assumption that $P$ is a mixture of $K$ nonparametric components that put most of their probability mass on sufficiently well separated subsets of the data space, we show that the components of the canonical functional resulting from the MLE of FMMs of ESDs correspond to the well separated mixture components of $P$. This can be interpreted as follows. If $P$ is a distribution generating well enough separated clusters of a very flexible kind, particularly not necessarily corresponding to the ESDs assumed to be the components of the FMM estimated by the MLE, the MLE will anyway for large enough data recover these clusters. We are not aware of any result of this kind in the literature regarding nonparametric consistency based on canonical functionals. 
\end{enumerate}
The paper is organized as follows.
In \Cref{sec:models_methods}, we define the class of ESDs and the corresponding FMMs. We review the connection between model-based clustering and FMMs.
We also define the MLE and give an account of the issue of unbounded likelihood. Finally we introduce and motivate the regularization of the ML optimization problem based on the ERC.
\Cref{sec:finite-samples} is devoted to the finite sample analysis of the ML procedure.
In \Cref{sec:consistency-1}, we establish the MLE's existence and consistency, assuming that $P$ is some general nonparametric distribution. 
\Cref{sec:consistency-2} investigates the canonical functional corresponding to the MLE for FMMS of ESDs for distributions $P$ that can be written as mixtures of sufficiently well separated nonparametric components.
In \Cref{sec:simulation} we present numerical experiments that illustrate the results.
\Cref{sec:conclusions} concludes the paper.


\section{Finite mixture modeling with elliptically symmetric distributions}
\label{sec:models_methods}
Elliptically symmetric distributions can be obtained applying an affine transformation to a spherical distribution. Let the random vector $Y \in \R^p$ have a spherical distribution \REV{with $\E{Y}={\bf 0}$ and $\Var{Y}\propto {\bf I}_p$, where ${\bf I}_p$ is the  $ p\times p$ identity matrix}. 
With fixed $\bmu \in \R^p$, and $\bOmega \in \R^{p{\times}p}$, the random vector $X = \bmu + \bOmega Y$ is said to have an elliptically-symmetric distribution denoted by $X \sim \ESD(\bmu, \bSigma)$, where $\bSigma \coloneqq  \bOmega \tr{\bOmega}$.
With a fixed distribution of the generating $Y$, $\ESD(\bmu, \bSigma)$ forms a symmetric location–scatter class of models with location parameter $\bmu$ and a symmetric positive semi-definite scatter matrix parameter $\bSigma$.
In particular  $\bmu = \E{X}$ and $\bSigma \propto \Var{X}$.
If $Y$ is absolutely continuous and $\bSigma$ is non-singular, than  $X$ is absolutely continuous with density function 
\begin{equation}\label{eq:esd}
f(\bx; \bmu, \bSigma)  \coloneqq  %
\det(\bSigma)^{-\frac{1}{2}} g\left(\tr{(\bx - \bmu)} \inv{\bSigma} (\bx - \bmu)\right),
\end{equation}
where $g:\Rnneg \to \Rnneg$ is the so-called radial or density generator
function.
Here we assume that $g(\cdot)$ is  monotonically non-increasing, which  implies that $f(\cdot)$ is unimodal.
\REV{Assume that $g(\cdot)$ is bounded, and}  $\lim_{t \to +\infty} g(t) = 0$ is required so that $f(\cdot)$ is a
proper density. 

The choice of $g(\cdot)$ (or, equivalently, the distribution of $Y$) defines a specific family of ESDs. A number of popular models can be obtained in this way, for example multivariate Gaussian distributions, multivariate Student-t distributions and multivariate logistic distributions.
Some families of ESDs have parameters in addition to $\bmu$ and $\bSigma$, \eg the degrees of freedom of the Student-t. We do not involve such additional parameters, so the results given here apply to MLEs for a family of multivariate Student-t distribution with fixed degrees of freedom, and generally to any location-scale family as defined in \cref{eq:esd} with $g$ fulfilling the conditions above.  

\cite{kelker70} originally introduced elliptical distributions to generalize the multivariate Gaussian model. ESDs have elliptically shaped density contours. They can model joint linear dependence between the $p$ components of a random vector and possibly heavy tails. For these reasons, ESDs occur frequently in the theory of statistics and applications \citep[for a comprehensive overview see][and references therein]{Paindaveine2014}. Finite mixture densities based on \cref{eq:esd} are the convex combinations
\begin{equation}\label{eq:psi}
\psi(\bx; \btheta)\coloneqq  \sum_{k=1}^K \pi_k f(\bx; \bmu_k, \bSigma_k),
\end{equation}
where $\pi_k \in [0,1]$ so that $\sum_{k=1}^K \pi_k = 1$, 
and $\btheta$ is a parameter vector collecting all elements of $\pi_k, \bmu_k, \bSigma_k,\; k=1,  \allowbreak  2, \allowbreak \ldots,  \allowbreak  K$. Given a family of ESDs, an FMM is defined by the set of distributions for all possible choices of $\btheta$.
FMMs are popular tools for modeling multimodality, skewness and heterogenous populations, and for performing several supervised and unsupervised tasks such as semi-parametric density estimation, classification and clustering \citep{mp00-book,fcr19-handbook}.

\subsection{Clustering and classification}
\label{sec:mix-clustering-classification}
FMMs used for 
model-based clustering and classification normally identify the classes or clusters with the
mixture components. FMMs can be equivalently formulated involving 
component membership indicators for the observations.

Assume that there are $K$ mixture components. Let $X_1,\ldots,X_n$ model a sequence of i.i.d. observations. These are accompanied by a sequence of 0-1 valued i.i.d. random variables $\zeta_i = \tr{(Z_{i1}, Z_{i2},\ldots,Z_{iK})},\ i=1,\ldots,n,$ with $\sum_{i=1}^K Z_{ik}=1$, i.e., for given $i$ one of the $Z_{ik}$ is 1, all the others are 0. Assume that $\zeta_i$ has a categorical distribution with $\prob{Z_{ik} = 1} = \E{Z_{ik}=1} = \pi_k$ for $k={1,  \allowbreak  2, \allowbreak \ldots,  \allowbreak  K}$.
Let \cref{eq:esd} be the density function of the conditional distribution of $X_i \given Z_{ik}=1$, then \cref{eq:psi} is the density function of the unconditional distribution of $X_i$. $Z_{ik}$ is the component (cluster) membership indicator for observation $i$ and the $k$th mixture component. \cref{eq:psi} can be seen as a model generating an expected fraction $\pi_k$ of points  from its $k$-th component $f(\cdot; \bmu_k, \bSigma_k)$.
If the $K$ mixture components are reasonably separated, sampling from \cref{eq:psi} will generate  clustered regions of data points. Cluster analysis is unsupervised, i.e., $\zeta_i,\ldots,\zeta_n$ are unobserved.

$\zeta_1,\ldots,\zeta_n$ represent a partition 
$\mathcal{G}_K \coloneqq  \set{G_k;\; k={1,  \allowbreak  2, \allowbreak \ldots,  \allowbreak  K}  }$ of $\{1,\ldots,n\}$, where 
$Z_{ik} \coloneqq  \ind\set{i \in G_k}$, $\ind\set{\cdot}$ being the usual indicator function, meaning that the mixture component (cluster) having generated $X_i$ is $k$. 

In cluster analysis one wants to assign an object to one of the $K$ groups of $\mathcal{G}_K$, which for $i=1,\ldots,n$ amounts to predicting $\zeta_i$ from $X_i$.

This can be done based on the posterior probability 
\begin{equation}\label{eq:tau}
\tau_k(\bx; \btheta) \coloneqq  \prob{Z_{ik} = 1 \given X_i = \bx} = %
\frac{\pi_k f(\bx; \bmu_k, \bSigma_k)}{\psi(\bx; \btheta)}.
\end{equation}
Given a parameter vector $\btheta$,
an  object $\bx$ can be assigned to the cluster $\cl(\bx,\btheta) \in \{1,  \allowbreak  2, \allowbreak \ldots,  \allowbreak  K\}$ according to the ``\emph{maximum a posteriori probability}'' (MAP) rule:
\begin{equation}\label{eq:map}
\cl(\bx,\btheta) \coloneqq  \argmax_{k \in \{1,  \allowbreak  2, \allowbreak \ldots,  \allowbreak  K\}}  \; \tau_k(\bx; \btheta). 
\end{equation}
The MAP rule implements the optimal Bayes classifier, \ie the assignment that  minimizes the expected 0-1 loss, also known as misclassification rate.
Since $\btheta$  is not known, it has to be estimated from the data, and \cref{eq:map} is then computed based on the estimator $\hat \btheta$.

\begin{remark}\label{rmk:optimality-map}
The optimality of the MAP rule for the clustering problem requires that (\emph{i}) data are generated from \cref{eq:psi}, (\emph{ii}) ``true'' clusters are defined in terms of  $Z_{ik} = \ind\set{i \in G_k}$, that is, $f(\cdot, \bmu_k, \bSigma_k)$ is the ``true''  underlying $k$-th class-conditional density, and (\emph{iii}) the posterior ratios in \cref{eq:tau} are computed at the ``true'' generating parameter $\btheta$ (the optimality would hold asymptotically with a consistent estimator for the ``true'' $\btheta$). In practice, these conditions are arguably never fulfilled ((ii) is a matter of definition, but relies on (i)). 

Applying this approach to data generated from a general distribution $P$ that does not necessarily have a density of type \cref{eq:psi} means that 
the method imposes a partition on $P$ that is governed by ``cluster prototype densities'' of the form  \cref{eq:esd}. Here we investigate what happens then, at least asymptotically.
\end{remark}

\subsection{Maximum likelihood estimator}
\label{sec:ml-intro}

The unknown mixture parameter vector $\btheta$ is typically estimated  by ML. Often, computations are performed based on the Expectation-Maximization (EM) algorithm \citep{mp00-book}.
Let $\mathbb{X}_n=\set{\bx_i,\; i={1,  \allowbreak  2, \allowbreak \ldots,  \allowbreak  n}   }$ be the observed sample. 
Define the sample likelihood and log-likelihood functions
\begin{equation}\label{eq:sample_lik}
\mathcal{L}_n(\btheta) := \prod_{i=1}^n \psi(\bx_i; \btheta), 
\quad 
\ell_n(\btheta) := \frac{1}{n} \sum_{i=1}^n \log(\psi(\bx_i; \btheta)). 
\end{equation}
Finding the maximum of  \cref{eq:sample_lik}, the sample MLE, is not straightforward.
In fact, \cref{eq:sample_lik} does not have a maximum. \cite{kw56}  discovered the unboundedness of $\ell_n$ in the context of univariate Gaussian FMM. The issue extends to many classes of FMMs, including those studied here. To see why, let us fix additional notations also used throughout the rest of the paper. Let $\norm{\cdot}$ be the Euclidean norm.
If the parameter vectors $\btheta$ come with indexes and accents, its members $\pi_k,\bmu_k,\bSigma_k$ have the same indexes and accents, \ie, $\tilde{\btheta}_m$ contains  $(\tilde{\pi}_{mk},\tilde{\bmu}_{mk},\tilde{\bSigma}_{mk})$ for all $m \in \N$.
For given $\btheta$, let  $\bkappa_k=(\bmu_k,\bSigma_k)$ (indexes and accents are applied to $\bkappa$ as above). 
Scalar parameters contained in $\btheta$'s sub-vectors are indexed so that the first two subscripts always denote the mixture component and the dimension in the feature space respectively, \eg   $\mu_{k,j} \in \R$ is the $j$-th coordinate of $\bmu_k$.

Consider a sequence $(\btheta_m)_{m\in\N}$ where $\bmu_{1,m}= \bx_1 \in \mathbb{X}_n$ and the smallest eigenvalue of $\bSigma_{1,m}$ converges to 0 as $m \conv +\infty$, then $\ell_n(\btheta_m) \conv +\infty$.
The likelihood degeneracy is caused by the fact that the density peak of $f(\cdot)$ is controlled by the smallest eigenvalue of the scatter matrix. 
In fact, \cref{eq:esd} can be parameterized in terms of the eigenvalue decomposition of $\bSigma$. That is
\begin{equation}\label{eq:esd_eigen}
f(\bx; \bmu, \bSigma) = %
\left( \prod_{j=1}^p \lambda_j(\bSigma) \right)^{-\frac{1}{2}} %
g\left( \sum_{j=1}^{p} \lambda_j(\bSigma)^{-1} \tr{(\bx-\bmu)}V_j(\bSigma)\tr{V_j(\bSigma)} (\bx-\bmu)\right), 
\end{equation}
where $\lambda_j(\bSigma)$ is the $j$-th eigenvalue of $\bSigma$, and $V_j(\bSigma)$ is its corresponding normalized eigenvector, \ie $\norm{V_j(\bSigma)} = 1$ for all $j= {1,  \allowbreak  2, \allowbreak \ldots,  \allowbreak  p}$.  Define  $\lmin^*(\bSigma)=\min\set{\lambda_{j}(\bSigma);\;  j={1,  \allowbreak  2, \allowbreak \ldots,  \allowbreak  p}}$, and  $\lmax^*(\bSigma) = \max\set{\lambda_{j}(\bSigma);\;  j={1,  \allowbreak  2, \allowbreak \ldots,  \allowbreak  p}}$.
The density $f(\cdot)$ can be bounded as
\begin{equation}\label{eq:esd_dens_upper_bound}
f(\bx; \bmu, \bSigma) %
\; \leq \; %
(\lmin^*(\bSigma))^{-\frac{p}{2}} g\left(\lmax^*(\bSigma)^{-1} \norm{\bx - \bmu}^2  \right)%
\; \leq \; %
g(0) (\lmin^*(\bSigma))^{-\frac{p}{2}}.
\end{equation}
Furthermore,
$f(\cdot) \in O(\lmin^{*}(\bSigma)^{-\frac{p}{2}})$, meaning that that $f(\bmu; \bmu, \bSigma) \conv +\infty$ as $\lmin^*(\bSigma) \searrow 0$ unless also $\lmax^*(\bSigma) \searrow 0$.
The ML problem can therefore not be solved in plain uncostrained form, and  
requires either constraints on the parameter space or a penalty.
Let $\Lambda(\btheta) = \{\lambda_{j}(\bSigma_k),\; \allowbreak j=1, \allowbreak \ldots, \allowbreak p,\; \allowbreak k=1, \allowbreak \ldots, \allowbreak k \}$, $\lmin(\btheta) = \min\{\Lambda(\btheta)\},\ \lmax(\btheta)=\max\{\Lambda(\btheta)\}$.

A possible constraint to define a proper ML problem
is the  eigenratio constraint (ERC) 
\begin{equation}\label{eq:gamma}
\frac{\lmax(\btheta)}{\lmin(\btheta)} \leq \gamma,
\end{equation}
with a constant $\gamma \ge 1$. $\gamma=1$ constrains all component scatter matrices to be spherical and equal.
Increasing $\gamma$ continuously allows for more flexible scatter shapes.
The ERC has been introduced by \cite{dennis81} and \cite{hathaway85} for the Gaussian case and brought back to the attention of the literature by \cite{ingrassia04}.
For a recent review see \cite{gggi18}.
The ERC captures the discrepancy between scatter shapes across components. Therefore, at least in clustering applications, these constraints have a data-analytic interpretation.
Another approach to solving the MLE existence problem is to rely on penalized ML methods. \cite{cri03} treated the case of the univariate Gaussian mixture. A recent comprehensive review is found in \cite{chen17}. 

In the following sections  we study the sequence of constrained MLEs 
\begin{equation}\label{eq:mle_sample}
\btheta_n \in  \argmax_{\btheta \in \tilde{\Theta}_K}\, \ell_n(\btheta),
\end{equation}
where the constrained parameter space is
\begin{equation}\label{eq:costrained-parameter-space}
\tilde{\Theta}_K:= \left\{
\btheta:  \; \; %
\pi_k \geq 0 \; \forall k\geq 1, \; %
\sum_{k=1}^K \pi_k =1; \; %
\frac{\lmax(\btheta)}{\lmin(\btheta)} \leq \gamma
\right\}.
\end{equation}
The ERC implies that $\lmax(\btheta) \in O(\lmin(\btheta))$ for all $\btheta \in \tilde{\Theta}_K$.
For a sequence $(\btheta_m)_{m \in \N}$ such that $\btheta_m \in \tilde{\Theta}_K$ and $\lambda_j(\bSigma_{k,m}) \searrow 0$  as $m \conv +\infty$ for some $j \in \{1,  \allowbreak  2, \allowbreak \ldots,  \allowbreak  p\}$ and $k \in \{1,  \allowbreak  2, \allowbreak \ldots,  \allowbreak  K\}$ the ERC enforces   $\lmax(\btheta_m) \searrow 0$. 
The ERC comes with two major technical challenges: (i) the parameter space $\tilde{\Theta}_K$ is not compact; (ii) the resulting ML problem cannot be solved by the  standard constrained optimization methods. 
\REV{
Note that, because of the ERC constraint, the estimator in \eqref{eq:mle_sample} is not affine equivariant.
However, affine equivariance can be made to hold if data are sphered first, see the \ref{sec:appendix-affine}.}
Other proposals of fully affine equivariant constraints exist in the literature 
\citep[see][]{gr09,ritter14-book}; however, there are no algorithms for their exact implementation.


\section{Finite sample existence}
\label{sec:finite-samples}

In this section we give precise non-asymptotic conditions involving $n,p,K$ and $g(\cdot)$ that guarantee the existence of $\btheta_n$ given the input data set $\mathbb{X}_n$.
The existence of the sequence $(\btheta_n)_{n \in \N}$  is the prerequisite for the study of its asymptotic behavior investigated in \Cref{sec:consistency-1} and \ref{sec:consistency-2}. 

\begin{assumption}\label{asm:fs}
For fixed $p,n,K, \mathbb{X}_n$:
\vspace{-\parskip}
\begin{enumerate}[label=(\alph*), ref=(\alph*)]

\item\label{asm:fs-n}
$n>K$, and $\mathbb{X}_n$ contains at least $K+1$ distinct points;

\item\label{asm:fs-og}
for all $\beta \in \Rpos$ and $\alpha \in \{0,1,\ldots, K\}$, 
$g(\beta y^{-1}  )^{n-\alpha} \in o(y^{\frac{p}{2}n})$ as $y \searrow 0$. 
\end{enumerate}
\end{assumption}

\Cref{asm:fs}-\ref{asm:fs-og} plays the central role in dealing with the sample likelihood function's unboundedness. 
It guarantees that for $\btheta \in \tilde{\Theta}_K$, such that  $\lmin(\btheta) \searrow 0$, the densities in the product $\mathcal{L}_n(\btheta)$ vanish sufficiently fast at all data points $\bx_i \in \mathbb{X}_n$ that do not coincide with mixture components' centers  $\bmu_k$ for all $k \in \{1,  \allowbreak  2, \allowbreak \ldots,  \allowbreak  K\}.$
The following Lemma is the key result to obtain the compactification of the parameter space.

\begin{lemma}\label{lem:fse_eigvals}
Let \Cref{asm:fs} hold. Consider a sequence $(\btheta_m)_{m \in \N} \in \tilde{\Theta}_K$  such that $\lambda_{j}(\bSigma_{k,m}) \searrow 0$ as $m \conv \infty$ for some $k \in \{1,  \allowbreak  2, \allowbreak \ldots,  \allowbreak  K\}$ and $j \in \{1,  \allowbreak  2, \allowbreak \ldots,  \allowbreak  p\}$.  Then, $\sup_{\tilde{\Theta}_K} \ell_n(\btheta_m) \conv  -\infty$  as $m \conv \infty$.
\end{lemma}

\begin{proof} %
First, rearrange $\mathcal{L}_n(\cdot)$. Consider a vector of indexes  $\bk := \tr{(k_1,k_2\ldots,k_n)}$ where $k_r \in \{1,  \allowbreak  2, \allowbreak \ldots,  \allowbreak  K\}$ for  $r \in \set{1,  \allowbreak  2, \allowbreak \ldots,  \allowbreak  n}$. For a fixed such $\bk$ define the products
\begin{equation}\label{eq:def_p_f_k}
{p}_{\bk}^f(\bx_1, \ldots, \bx_n; \btheta):= \prod_{r=1}^n f(\bx_r; \bmu_{k_r}, \bSigma_{k_r}), 
\quad \text{and} \quad 
p^{\pi}_{\bk}(\btheta):= \prod_{r=1}^{n} \pi_{k_r}.
\end{equation}

There exists $K^n$ of such possible vectors of indexes $\bk$, say $\bk_1, \bk_2, \ldots, \bk_{K^n}$. Based on these vectors it is possible to write the sample likelihood function as a mixture of $K^n$ components like \cref{eq:def_p_f_k}, that is 
\begin{equation}\label{eq:Ln_mix_products}
\mathcal{L}_n(\btheta) = \sum_{h=1}^{K^n}  p^{\pi}_{\bk_h}(\btheta)\, p^{f}_{\bk_h}(\bx_1, \ldots, \bx_n; \btheta).
\end{equation}
Since $\btheta_m \in \tilde{\Theta}_K$, $\lambda_{j}(\bSigma_{km}) \searrow 0$ implies that $\lmax(\btheta_m) \searrow 0$ as $m \conv +\infty$ 
Assume w.l.o.g. (otherwise consider a suitable subsequence) that for sufficiently large $m$,  the sequence $(\btheta_m)_{m \in \N}$ is such that, for $j \in \{1,  \allowbreak  2, \allowbreak \ldots,  \allowbreak  p\},$ and $k \in \{1,  \allowbreak  2, \allowbreak \ldots,  \allowbreak  K\}$, the centrality parameter $\mu_{kjm}$  either converges, or it leaves any compact set. %
If the limits of $(\bmu_{km})_{m \in N}$ belong to $\mathbb{X}_n$, there are at most $K$ of them.
Therefore, according to \Cref{asm:fs}-\ref{asm:fs-n}, there exists $i \in \{1,  \allowbreak  2, \allowbreak \ldots,  \allowbreak  n\}$ and $\eps>0$ such that $\bx_i\in \mathbb{X}_n{\setminus}\mathbb{X}'$ and  $\norm{\bx_i - \bmu_{km}} > \eps$ for large enough $m$. %
Consider a factorization such as the one in \cref{eq:def_p_f_k} for some  vector of indexes $\bk$.  Define $\mathbb{X}'_{\bk} \subseteq \mathbb{X}'$ where 
\[
\mathbb{X}'_{\bk} :=\set{\bx_r \in \mathbb{X}': \lim_{m \to \infty} \bmu_{k_r,m} = \bx_r,\; \text{for any}\; k_r\in\{k_1,k_2, \ldots, k_n\}}.
\]
Let $\#(\mathbb{X}'_{\bk}) = q_{\bk}\le K$. Therefore, $p_{\bk}^{f}(\cdot)$ in  \cref{eq:def_p_f_k} can be factorized as
\begin{equation}~\label{eq:factor_p_f_k}
{p}_{\bk}^f(\bx_1, \ldots, \bx_n; \btheta_m)  =    %
   \prod_{r:\,\bx_r \in \mathbb{X}'_{\bk}}  %
      f(\bx_r; \bmu_{{k_r}m}, \bSigma_{{k_r}m}) 
   \prod_{r:\,\bx_r \in \mathbb{X}{\setminus}\mathbb{X}'_{\bk}}
      f(\bx_r; \bmu_{{k_r}m}, \bSigma_{{k_r}m}).
\end{equation}
As $m \conv +\infty$,   $f(\bx_r; \bmu_{{k_r}m}, \bSigma_{{k_r}m}) \in O(\lmin(\btheta_m)^{-\frac{p}{2}})$ for all $r$ such that $\bx_r \in \mathbb{X}'_{\bk}$, therefore
\[
\sup_{\tilde{\Theta}_K} \; \prod_{r:\,\bx_r \in \mathbb{X}'_{\bk}}  %
f(\bx_r; \bmu_{{k_r}m}, \bSigma_{{k_r}m}) %
\in O\left( \lmin(\btheta_m)^{-\frac{p}{2}q_{\bk}}   \right).
\]
On the other hand, for all $r$ such that $\bx_r \in \mathbb{X}{\setminus}\mathbb{X}'_{\bk}$ and a positive constant $c_{k_r}$
\[
f(\bx_r; \bmu_{{k_r}m}, \bSigma_{{k_r}m}) %
\leq %
O(\lmin(\btheta_m))^{-\frac{p}{2}} \; g(\lmin(\btheta_m)^{-1} c_{k_r} ).
\]
\Cref{asm:fs}-\ref{asm:fs-og} and  $q_{\bk} \leq K$ ensure that  $g(\lmin(\btheta_m)^{-1} c_{k_r} )^{n-q_{\bk}} \in o(\lmin(\btheta_m)^{\frac{p}{2}n})$, therefore 
\[
\sup_{\tilde{\Theta}_K} \; %
\prod_{r:\,\bx_r \in \mathbb{X}{\setminus}\mathbb{X}'_{\bk}} f(\bx_r; \bmu_{{k_r}m}, \bSigma_{{k_r}m}) %
\in %
O \left(  \lmin(\btheta_m)^{-\frac{p}{2}{(n-q_{\bk})}} \right) %
\; o\left( \lmin(\btheta_m)^{\frac{p}{2}n} \right).
\]
The latter implies that
\[
\sup_{\tilde{\Theta}_K} \; {p}_{\bk}^f(\bx_1, \ldots, \bx_n; \btheta_m) %
\in %
O\left( \lmin(\btheta_m)^{-\frac{p}{2}q_{\bk}}   \right) %
O \left(  \lmin(\btheta_m)^{-\frac{p}{2}{(n-q_{\bk})}} \right) %
o\left( \lmin(\btheta_m)^{\frac{p}{2}n} \right)
\]
That is
\[
\sup_{\tilde{\Theta}_K} \; {p}_{\bk}^f(\bx_1, \ldots, \bx_n; \btheta_m) %
\in o(1).
\]
We conclude that $\sup_{\tilde{\Theta}_K} {p}_{\bk}^f(\bx_1, \ldots, \bx_n; \btheta_m) \conv 0$ for large enough  $m$, whatever the vector of indexes $\bk$. The latter implies that all factors of $\mathcal{L}_n(\btheta_m)$ in   \eqref{eq:Ln_mix_products} vanish and therefore $\sup_{\tilde{\Theta}_K} \ell_n(\btheta_m) \conv  -\infty$.
\end{proof}

\begin{theorem}[finite sample existence]\label{th:fs_existence}
Under \Cref{asm:fs}, $\btheta_n$ exists.
\end{theorem}
\begin{proof}
First note that $\tilde{\Theta}_K$ is not empty because for any $\gamma \geq 1$  and any choice  $\bSigma_1, \bSigma_2, \ldots \bSigma_K$ such that the ERC is not fulfilled, for all $k=1,2, \ldots, K$ one can always replace $\bSigma_k$ with 
$\dot \bSigma_k = V(\bSigma_k) \dot{\Lambda}(\bSigma_k) \tr{V(\bSigma_k)}$, where $V(\bSigma_k) \in \R^{p\times p}$ is an orthogonal matrix whose columns are eigenvectors of $\bSigma_k$, and  $\dot{\Lambda}(\bSigma_k) \in \R^{p \times p}$ is a diagonal matrix where
$\dot{\Lambda}(\bSigma_k)[j,j] = \min\{\lambda_j(\bSigma_k), \gamma \lmin(\btheta)\}$. By construction any such choice $\dot \bSigma_1, \dot \bSigma_2, \ldots \dot \bSigma_K$ will always satisfy the eigen-ratio constraint.
With the following step we show that there exists a compact set $T_K \subset \tilde{\Theta}_K$ such that $\sup_{\btheta \in T_K} \ell_n(\btheta) = \sup_{\btheta \in \tilde{\Theta}_K} \ell_n(\btheta)$.

{Step (a):} %
Consider $\btheta$ such that  $\pi_1=1$, $\bmu_1= \bx_1$, $\bSigma_j=\bI_p$ for all $k\in\{1,2,\ldots,K\}$, arbitrary $\bmu_k$ and $\pi_k=0$ for all $k \neq 1$.  For such a choice of $\btheta$, $l_n(\btheta) = \inv{n} \sum_{i=1}^n \log f(\bx_i; \bx_1,\bI_p) > -\infty$, thus $\sup_{\btheta \in \tilde{\Theta}_K} l_n(\btheta) > -\infty$.\\
{Step (b):} %
Consider a sequence $(\dot \btheta_m )_{m\in N}$.  Lemma \ref{lem:fse_eigvals} rules out the possibility that $\dot \lambda_{kjm} \searrow 0$ for some index  $k \in \{1,  \allowbreak  2, \allowbreak \ldots,  \allowbreak  K\}$ and $j \in \{1,  \allowbreak  2, \allowbreak \ldots,  \allowbreak  p\}$, because this would imply that $\sup_{\tilde{\Theta}_K} \ell_n(\dot \btheta_m) \conv  -\infty$.
Using \cref{eq:esd_dens_upper_bound},
\begin{equation}\label{eq:sample-loglik-upper-bound}
\ell_n(\btheta) =    %
\frac{1}{n} \sum_{i=1}^n \log \left( \sum_{k=1}^K \pi_k f(\bx_i; \bmu_k, \bSigma_k)  \right) %
\leq K \log \left(2\pi\lambda_{\min}(\btheta)  \right)^{-\frac{p}{2}}. 
\end{equation}
Assume that $\dot \btheta_m  \conv \dot \btheta$ where $ \dot \btheta \in \tilde{\Theta}_K$ and  $\dot{\lambda}_{kjm} \conv +\infty$ for some indexes  $k \in \{1,  \allowbreak  2, \allowbreak \ldots,  \allowbreak  K\}$ and $j \in \{1,  \allowbreak  2, \allowbreak \ldots,  \allowbreak  p\}$.
Because of the ERC, $\lambda_{\min}(\dot \btheta_m) \conv +\infty$, and \cref{eq:sample-loglik-upper-bound} implies that $\sup_{\dot \btheta_m} \ell_n(\dot \btheta_m) \conv -\infty$. Therefore, the case that $\dot{\lambda}_{kjm} \conv +\infty$  is also ruled out.\\
{Step (c):} %
now suppose that  $\norm{\dot \bmu_{km}} \conv +\infty$ for some $k \in \{1,  \allowbreak  2, \allowbreak \ldots,  \allowbreak  K\}$. W.l.o.g take $k=1$. Choose an alternative sequence $(\ddot \btheta_m)_{m \in \N}$ that is equal to $(\dot \btheta_m)_{m \in \N}$ except now $\ddot \bmu_{1m} = \bm{0}$ for all $m \in \N$. Note that $f(\bx_i;  \dot \bmu_{1m}, \dot \bSigma_1) \conv 0$ for all $i \in \{1,  \allowbreak  2, \allowbreak \ldots,  \allowbreak  n\}$, which implies that  $\psi(\bx_i, \dot \btheta_m) < \psi(\bx_i,\ddot \btheta_m)$ for large enough $m$ for all $i \in \{1,  \allowbreak  2, \allowbreak \ldots,  \allowbreak  n\}$. Hence, \REV{$\ell_n(\dot \btheta_m) < \ell_n(\ddot \btheta_m)$} for large $m$. 

Steps (a)--(c) plus the fact that $\pi_k \in [0,1]$ for all $k \in \{1,  \allowbreak  2, \allowbreak \ldots,  \allowbreak  K\}$ imply that the vector $\btheta$ maximizing $\ell_n(\cdot)$ must lie in a compact subset $T_n \subset \tilde{\Theta}_K$, and the continuity of $\ell_n(\cdot))$ guarantees existence of $\btheta_n$. 
\end{proof}

\begin{remark}[Finite sample existence for specific models]\label{rmk:popular-finite-sample-cases}
The key \Cref{asm:fs}-\ref{asm:fs-og} holds depending on the density generator function $g(\cdot)$.
For the Gaussian distribution, $g(t) = c_p \exp(-t/2)$, with $c_p$ being a  constant dependent on $p$.
It is easy to see that  the condition is fulfilled for all $n>K$.
For the Student-t distribution with $\nu$ degrees of freedom,
$g(t) = c_{p,\nu} (1 + t/\nu)^{-(\nu+p)/2}$ where $c_{p,\nu}$ is a constant that depends on both $p$ and $\nu$.
In the latter case \Cref{asm:fs}-\ref{asm:fs-og} holds if $n > K ( 1 + p/\nu)$, requiring more data points for larger $p$ and smaller $\nu$. Not surprisingly,  if $\nu \conv +\infty$, the condition is fulfilled with $n>K$ as for the Gaussian case.
\Cref{asm:fs}-\ref{asm:fs-og} can be easily checked for other ESDs as well.
\end{remark}

\begin{remark}[Computing]\label{rmk:em-algorithm}
The previous statement is also relevant to ensure that applying computing algorithms searching for $\btheta_n$ on a given input data set makes sense.
In the FMM context, the usual approach to compute the MLE is to run the Expectation-Maximization (EM) algorithm of \cite{delaru77} or some of its variants \citep[see][]{mk97-book}.
\cite{ggmm15} and \cite{ch17-jmlr} developed EM-type algorithms implementing the ERC for the case of a Gaussian FMM;  \cite{ch17-jmlr} also proved convergence results.
The general structure of the EM applied to the FMM will apply to finite mixtures of ESDs; however, the generality of $g(\cdot)$ in \cref{eq:esd} does not allow  to write down the M-step explicitly.
Furthermore, the implementation of the ERC for the Gaussian case can be extended to those ESD that have a Gaussian representation. 
For instance, for mixtures of Student-t distributions, one can take the EM algorithm presented in Chapter 7 of \cite{pm00} and add the ERC via a conditional M--step similar to the  CM1-step of Algorithm 2 in \cite{ch17-jmlr}. 
In general, an appropriate constrained M-step needs to be developed depending on the specific $g(\cdot)$, but this is outside the scope of this paper.
\end{remark}


\section{Asymptotic analysis for general $P$}
\label{sec:consistency-1}

Now consider a general $P$ as generator of $p$-dimensional data to be analysed by an MLE derived from an FMM of ESDs. Treating $K$ in the definition of $\psi$ as fixed,
define the following log-likelihood-type functionals: 
\[
L(\btheta,P) \coloneqq   \int\log\psi(\bx,\btheta)dP(\bx),
\]
and 
\[
L_K \coloneqq L_K(P) \coloneqq \sup_{\btheta\in\tilde{\Theta}_K} L(\btheta),%
\qquad 
\btheta^\star \coloneqq \btheta^\star(P) \coloneqq \argmax_{\btheta \in \tilde{\Theta}} L(\btheta,P).
\]
$\btheta^\star(P)$ is not normally unique (in particular, there is the ``label switching'' issue, i.e., the numbering of the mixture components is arbitrary), and is taken to be any of the maximising parameter vectors.  
Without demanding that $P$ is of type \cref{eq:psi}, the usual interpretation of these quantities is that $L(\btheta; P)$ is proportional, up to a term that only depends on $P$, to the Kullback-Liebler risk (KLR) of approximating the density of $P$ by $\psi(\cdot, \btheta)$.
Therefore,  $\btheta^\star(P)$  provides the best approximation to $P$ in terms of KLR  obtainable from a $K$-components ESD mixture model. The following analysis establishes the existence of $\btheta^\star(P)$ and the convergence of $(\btheta_n)_{n \in \N}$ to $\btheta^\star(P)$. The notation $\E[P]{\cdot}$ denotes the expectation with respect to the distribution $P$.

\begin{assumption}\label{asm:no-mass}
For every set $A=\set{\bx_{1},\bx_{2},\ldots,\bx_{K}} \subset \R^{p}$ with at most $K$ points: $P(A) < 1.$
\end{assumption}

\begin{assumption}\label{asm:bounded-ElogGx2}
$\E[P]{  \log g(\norm{X}^{2})   } > -\infty$. 
\end{assumption}

\begin{assumption}\label{asm:bounded-Elogf}
For fixed $\beta \in \Rpos$, $g(\beta y^{-1}  ) \in o(y)$ as $y \searrow 0$. 
\end{assumption}

\begin{assumption}\label{asm:monotone-LK}
$L_{K-1}(P) < L_{K}(P)$.
\end{assumption}

\Cref{asm:no-mass} is fulfilled  if $P$ does not concentrate its mass on $K$ points.
Note that this assumption is fulfilled by the empirical distribution $P_n$ for sufficiently large $n$ when $P$  is absolutely continuous.
\Cref{asm:bounded-ElogGx2} is needed to ensure that under $P$ it makes sense to maximize the log-likelihood function. When $f(\cdot)$ is the Gaussian density,  \Cref{asm:bounded-ElogGx2} is fulfilled if $\oE_P \norm{X}^2$ exists. 
\Cref{asm:bounded-Elogf} is the analog of \Cref{asm:fs}-\ref{asm:fs-og} and deals with the degeneration of the scatter matrices.
It implies that, far from the centers $\bmu_k$,  $f(\cdot)$ vanishes at a speed that is sufficiently fast compared to the speed at which it becomes unbounded in the center when $\lmin^*({\bSigma_k}) \searrow 0$.
\Cref{asm:bounded-Elogf} is fulfilled by the Gaussian model.
When $f(\cdot)$ is the Student-t distribution with $\nu$ degrees of freedom, it holds if $\nu + p > 2$.


Regarding \Cref{asm:monotone-LK}, note that generally $r<s\Rightarrow L_{r}(P) \le L_{s}(P)$, because any parameter $\btheta$ with $r$ mixture components can be reproduced with more mixture components setting some component proportions to 0. If $L_{K-1}(P) = L_{K}(P)$, then maxima of the likelihood exist with a component proportion of 0, and the corresponding location and scatter parameters can take any value. Particularly then the maxima of the likelihood cannot be forced  into any compact set. 


To prove the consistency \Cref{thm:mle_consistency}, we first establish the existence of the ML functional $\btheta^\star$ (or equivalently $\btheta^\star(P)\neq\emptyset$) in \Cref{thm:asy_existence}. The existence \Cref{thm:asy_existence} is obtained via the compactification of the parameter space based on the following Lemmas \ref{lem:LK-gt-inf}-\ref{lem:bounded-mu}.

\begin{lemma}\label{lem:LK-gt-inf}
Under \Cref{asm:bounded-ElogGx2},  for all $K\geq 1$, $L_{K}(P) > -\infty.$
\end{lemma}

\begin{proof} 
Choose $\btheta \in \tilde{\Theta}_K$ such that $\bmu_1=\bm{0}$, $\bSigma_1=\bI_p$, $\pi_1=1, \pi_2=\ldots=\pi_K=0$. Note that $\bSigma_1$ fulfills the ERC.
All remaining parameters of $\btheta$ are chosen arbitrarily. The statement is proven by using \Cref{asm:bounded-ElogGx2}:
\[
\sup_{\btheta \in \tilde{\Theta}_K}\int \log \psi(\bx, \btheta)dP(\bx) \geq
\int \log \pi_1 f(\bx; \bmu_1, \bSigma_1)dP(\bx) \geq 
\int \log g(\norm{\bx}^2)dP(\bx) > -\infty.
\]
\end{proof}

\begin{lemma}\label{lem:bounded-eig}
Under \Cref{asm:no-mass,asm:bounded-ElogGx2,asm:bounded-Elogf}  %
there exist  $\lambda^\star_{\min}>0,\ \lambda^\star_{\max}<\infty,\ \epsilon>0$,  so that
\begin{enumerate}[label=(\alph*), ref=(\alph*)]

\item\label{lem:bounded-eig-L}
$L(\btheta, P) \le L_K-\epsilon$ for every $\btheta$ with  $\lmin(\btheta)<\lmin^\star$ or $\lmax(\btheta) > \lmax^\star$,

\item\label{lem:bounded-eig-ln}
for iid samples $X_1,X_2,\ldots$ from $P$, for sequences $(\dot \btheta_n)_{n \in \N}$ with 
$\lmin(\dot \btheta_n)<\lmin^\star$ or $\lmax(\dot \btheta_n)>\lmax^\star$ for large enough $n$: 
$\ell_n(\dot \btheta_n) \le \ell_n(\btheta_n) - \epsilon$ $P$-a.s.
\end{enumerate}
\end{lemma}

\begin{proof} 
Consider a sequence $(\btheta_m)_{m\in\N}$ where $\btheta_m \in \tilde{\Theta}_K$ with $\lmax(\btheta_m) \conv +\infty$  for all  $m \in \N$.
The  ERC enforces $\lambda_{\min}(\btheta_m) \conv +\infty$ and therefore $\sup_{\bx} f(\bx, \bmu_{km},\bSigma_{km}) \searrow 0 $ for all $k=1,2,\ldots, K,$ and   $\sup_{\bx} \psi(\bx; \btheta_m) \searrow 0$. The dominated convergence theorem implies that $\E[P]{\psi(\bx; \btheta_m)} \searrow 0$, and therefore  $L(\btheta_m, P)  \leq  \log( \E[P]{\psi(\bx; \btheta_m)}) \conv -\infty$. $L(\btheta_m, P)\searrow -\infty$  according to Lemma \ref{lem:LK-gt-inf} makes it impossible that $L(\btheta_m, P)$ is close to $L_K(P)$ for $m$ large  enough. The latter proves the existence of  the upper bound $\lmax^\star < \infty$ as required in part (a).\\

Now consider a sequence $(\btheta_m)_{m\in\N}$ where $\btheta_m \in \tilde{\Theta}_K$ with  $\lmin(\btheta_m) \searrow 0$.
Define 
\[
A_{m,\epsilon} = \set{ \bx: \; \min_{1 \leq k \leq K} \|\bx-\bmu_{km}\|  > \epsilon }.
\]
\Cref{asm:no-mass} ensures that for $\delta>0$ there exists $\epsilon>0$ so that $P(A_{m,\epsilon}) \ge \delta $ for all $m\in\N$.
\begin{align*}
L(\btheta_m, P) \leq %
         &\int_{A_{m,\epsilon}} \log \psi(\bx; \btheta_m)dP(\bx) + \int_{A^c_{m,\epsilon}} \log \psi(\bx; \btheta_m)dP(\bx)\\
{}\leq   &\int_{A_{m,\epsilon}}  \log \left( K (2\pi)^{-\frac{p}{2}} \lmin(\btheta_m)^{-\frac{p}{2}} g(\gamma\lmin(\btheta_m)^{-1} \epsilon ^ 2) \right)dP(\bx) + \\
{}     &+ \int_{A^c_{m,\epsilon}}  \log \left(K(2\pi)^{-\frac{p}{2}}  \lmin(\btheta_m)^{-\frac{p}{2}} g(0)   \right)dP(\bx) \\
{}\leq   &{~}P(A_{m,\epsilon})\left(\log K - \frac{p}{2}\log(2\pi) - \frac{p}{2} \log \lmin(\btheta_m) + \log (g(\gamma \lmin(\btheta_m)^{-1} \epsilon^2) \right) + \\
{}      &+ P(A^c_{m,\epsilon}) \left( \log K - \frac{p}{2}\log(2\pi) - \frac{p}{2} \log \lmin(\btheta_m) + \log g(0)\right)\\
{}\leq   &{~} c_1 + c_2 \log \lmin(\btheta_m)^{-1} + c_3 \log g\left(  c_4 \lmin(\btheta_m)^{-1} \right)
\end{align*}
for positive constants $c_1, c_2,c_3$ and $c_4$ all independent of $\btheta_m$. 
The previous inequality can be rearranged as follows
\[
L(\btheta_m, P) \leq
\log g\left(  c_4 \lmin(\btheta_m)^{-1} \right)
\left( %
\frac{c_1}{\log g\left(  c_4  \lmin(\btheta_m)^{-1} \right)}  %
- c_2 \frac{\log \lmin(\btheta_m)}{\log g(c_4\lmin(\btheta_m)^{-1}) } %
+ c_3
\right).%
\]
By  \Cref{asm:bounded-Elogf}, as $m \conv +\infty$,  $g(c_4\lmin(\btheta_m)^{-1}) \searrow 0$ faster than $\lmin(\btheta_m)$ and   $L_K(P) \conv -\infty$, which contradicts $L_K(P)>-\infty$ established by Lemma \ref{lem:LK-gt-inf}.
Therefore, there exists a lower bound $\lmin^\star >0$ as stated in part (a) of the statement.

Regarding part (b),  let the sequence $(\dot \btheta_n)_{n\in\N}$ be  chosen as the sequence above taking $m=n \conv \infty$. Let $P_n$ be the empirical distribution.
The class of all $A_{n,\epsilon}$ is a subset of the class of  intersections of the complements of all closed balls, and therefore a Vapnik-Chervonenkis class  \citep[see][]{vw96-book}.
Glivenko-Cantelli enforces $P_n(A_{n,\epsilon}) - P(A_{n,\epsilon})\to 0$ a.s.
Note that  $L(\dot \btheta_n, P_n) = \ell_n(\dot \btheta_n)$, and $L_K(P_n) = \ell_n(\btheta_n)$.
For large enough $n$ Lemma \ref{lem:LK-gt-inf} holds with $P_n$ replacing $P$. Therefore part (b) of the statement  is proved by replacing $P$ with $P_n$ in the proof of part (a).
\end{proof}

\begin{lemma}\label{lem:bounded-mu}
Under \Crefrange{asm:no-mass}{asm:monotone-LK}, there is a compact set $T\subset\R^p$, $\epsilon_1 >0$ and $\epsilon_2 >0$ so that
\begin{enumerate}[label=(\alph*), ref=(\alph*)]
\item\label{lem:bounded-mu-L}
$L(\btheta, P) \leq L_K - \epsilon_1$ whenever $\bmu_k \notin T$ for some $k \in \{1,  \allowbreak  2, \allowbreak \ldots,  \allowbreak  K\}$;

\item\label{lem:sup-bounded}
$\sup_{\btheta \in \tilde{\Theta}_{K} : \bmu_1,\ldots,\bmu_K \in T} L(\btheta, P) = L_K(P) < +\infty$;

\item\label{lem:bounded-mu-ln}
for iid samples $X_1,X_2,\ldots$ from $P$, for sequences $(\dot \btheta_n)_{n \in \N}$ with $\dot \btheta_n \in \tilde{\Theta}_{K}$ and $\dot \bmu_{kn} \in T$ for all $k \in \{1,  \allowbreak  2, \allowbreak \ldots,  \allowbreak  K\}$  for large enough $n$:
$\sup_{\dot \btheta_n} \ell_n(\dot \btheta_n) \le \ell_n(\btheta_n) - \epsilon_2$ $P$-a.s
whenever $\dot \bmu_{kn} \notin T$ for some $k \in \{1,  \allowbreak  2, \allowbreak \ldots,  \allowbreak  K\}$,
and
$\sup_{\dot \btheta_n} \ell_n(\dot \btheta_n) = \ell_n(\btheta_n)$  $P$-a.s.
\end{enumerate}
\end{lemma}

\begin{proof} 
Consider a sequence $(\btheta_m)_{m\in\N}$ where $\btheta_m \in \tilde{\Theta}_K$ for all $m \in \N$.
Assume $\norm{\bmu_{km}} \conv  \infty$ for $k \in \{1,2,\ldots, r\}$, and $\bmu_{km}\in T$ for all $k>r$.
Note that $r=K$ would imply $L(\btheta_m, P) \conv -\infty$ for large enough $m$, therefore, we only consider the case when $r < K$.
Take a sequence $(\epsilon_m)_{m \in \N}$ and construct the sets 
\begin{equation*}
A_m:=\set{
x: \;  \forall k \in \set{1,\ldots,r} \;
f(\bx,\bmu_{km},\bSigma_{km})
\leq
\epsilon_m \sum_{j = r+1}^K \pi_{jm} f(\bx,\bmu_{jm},\bSigma_{jm})
}, 
\end{equation*}
where  $\epsilon_m\searrow 0$ slowly enough that $A_1 \supset A_2 \supset, \ldots$ and 
$P(A_m) \nearrow 1$.
Define another sequence $(\bar \btheta_{m})_{m \in \N}$ such that $\bar \btheta_{m} \in \tilde{\Theta}_{K-r}$ with $\bar \bmu_{km} = \bmu_{(k+r)m}$,  
$\bar \bSigma_{km} = \bSigma_{(k+r)m}$, and $\bar \pi_{km} = \pi_{(k+r)m}\inv{(\sum_{j = r+1}^K \pi_{jm})}$ for all $k \in \{r+1, \ldots, K\}$.
By construction $\bar \bmu_{km} \in T$ for all $k \in \{r+1, \ldots, K\}$.
Lemma \ref{lem:bounded-eig} implies that for all $\btheta \in \tilde{\Theta}_K$ that are sufficiently close to $L_L$ (that is $L(\btheta, P) > L_k - \epsilon$ according to Lemma \ref{lem:bounded-eig}) the mixture component densities are bounded above: 
$f(\bx; \bmu_k, \bSigma_k) \leq f_{\max} = (2\pi \lmin^{\star})^{-\frac{p}{2}}$
for all $k$ and some suitably defined $\lmin^{\star}>0$. The latter also implies that $L_K < +\infty$.  $\psi(\bx,\bar \btheta_{m}) \geq  \psi(x, \btheta_{m})$ for all  $\bx \in A_m$, therefore
\begin{align*}
L(\btheta_m,P)
   &= \int_{A_m}\log(\psi(\bx,\btheta_{m})) dP(\bx)
   + \int_{A_m^c}\log(\psi(\bx,\btheta_{m}))dP(\bx) \nonumber \\
   &\leq \int_{A_m}\log ((1+\epsilon_m)\psi(\bx, \bar \btheta_{m})) dP(\bx)
   + P(A_m^c) \log (f_{\max}). \\
   &\leq \int \ind_{A_m}(\bx) \log(\psi(\bx, \bar \btheta_{m})) dP(\bx)) + a_m
   \end{align*}
   with  $a_m \searrow 0$.  $\log(f(\bx, \bar \bmu_{km}, \bar \bSigma_{km})$ can be bounded by $\log(f_{\max})$, and by applying the dominated convergence theorem we obtain that $ \int \ind_{A_m}(\bx) \log(\psi(\bx, \bar \btheta_{m})) dP(\bx) \conv L(\bar \btheta_m, P)$. Therefore,  $L(\btheta_m, P) \leq  L(\bar \btheta_m, P) + o(1)$  for large enough $m$. Because of \Cref{asm:monotone-LK}, $L(\bar \btheta_{m},P) \leq L_{K-r} < L_K$ and therefore  $L(\btheta_m,P) < L_K$ for sufficiently large $m$, which proves part \ref{lem:bounded-mu-L} of the statement.

Observe that Lemma \ref{lem:bounded-eig} implies that for all $\btheta \in \tilde{\Theta}_K$ for which $L(\btheta)$ is close to $L_K$, each eigenvalue in  $\btheta$ is contained in $[\lmin^\star, \lmax^\star]$ with $0 < \lmin^\star \leq \lmax^\star < +\infty$, 
$\psi(\bx,\btheta) \leq f_{\max}$, and therefore $L_K  < +\infty$.
Consider   $(\btheta_m)_{m \in \N}$ so that  $L(\btheta_m) \conv L_k$ with $\bmu_{1m},\ldots,\bmu_{Km}\in T$, and  $0 < \lmin^\star \leq \lambda_{kjm} \leq \lmax^\star < +\infty$ for all  $k \in \{1,\ldots, K\}$,  $j \in \{1,\ldots,p\}$.
Because of compactness, there exists $\btheta_0$ such that $\btheta_m \to \btheta_0$ and, using Fatou's Lemma, $L_K= \lim_{m\to\infty} L(\btheta_m,P) \leq \E[P]{\limsup_m \psi(\bx,\btheta_m)} = L(\btheta_0) \leq L_K < +\infty$. The latter proves the existence of a maximum stated in part \ref{lem:sup-bounded}.

As for the proof of Lemma \ref{lem:bounded-eig}-\ref{lem:bounded-eig-ln}, part  \ref{lem:bounded-mu-ln} of the statement is shown by setting $n=m$, replacing the previous sequences $(\btheta_m)_{m\in\N}$ and $({\bar{\btheta}}_m)_{m\in\N}$ with $(\dot \btheta_n)_{n\in\N}$ and $(\dot{\bar{\btheta}}_n)_{n\in\N}$, and replacing $P$ with  the empirical distribution $P_n$.
The same steps as before can be applied taking into account the following observations.
Glivenko-Cantelli enforces $P_n(A_{n})-P(A_{n}) \conv  0$ a.s. 
because we can construct a sequence of closed balls $(B_n)_{n\in \N}$  so that $B_n\subseteq A_n$ and  $P(B_n) \conv  1$ a.s. 
The closed balls are a Vapnik-Chervonenkis class, and $P_n(A_{n}) \geq P_n(B_n)$.
Note that for all  $\btheta \in \tilde{\Theta}_K$ with $L(\btheta, P) > L_{K-1}$: $\ell_n(\btheta) \conv  L(\btheta, P)$ a.s. by the strong law of large numbers and therefore for large enough $n$: $\sup_{\btheta \in \tilde{\Theta}_K} \ell_n(\btheta) > L_{K-1}$ a.s.
Furthermore, if $\dot{\bar{\btheta}}_{n}$  is chosen optimally in a compact set as in the proof of part \ref{lem:sup-bounded},  $\ell_n(\cdot)$ converges uniformly over $\tilde{\Theta}_{K-r}$ (see Theorem 2 in \cite{Jennrich69}) and therefore 
$\limsup_{n\to\infty} \ell_n(\dot{\bar{\btheta}}_{n}) \leq L_{K-1} < L_K$.
\end{proof}

\begin{theorem}[Existence of the ML functional]\label{thm:asy_existence}
Under \Crefrange{asm:no-mass}{asm:monotone-LK} there is a compact subset  $T_{\tilde{\Theta}_K} \subset \tilde{\Theta}_K$ so that there exists $\btheta\in T_{\tilde{\Theta}_K}$ such that $-\infty < L(\btheta) = L_{K} < +\infty$, and  for all $\btheta \not\in T_{\tilde{\Theta}_K}$, $L(\btheta)$ is bounded away from $L_K$.
\end{theorem}

\begin{proof}
The statement is shown by putting together Lemmas \ref{lem:LK-gt-inf}--\ref{lem:bounded-mu}.
\end{proof}

\cite{hmg06} found sufficient conditions for the identifiability of certain classes of finite mixtures of ESDs.
If $P=P_0$, and $\psi(x; \btheta_0)$ is its density whose mixture components also fulfill the sufficient conditions of \cite{hmg06}, the previous theorem guarantees that the argmax functional $\btheta^{\star}(P)$ exists and is unique up to label switching. 
This is a consequence of Lemma 1 of \cite{wald49}.
But here we are interested in the more realistic case when $\psi$ is not necessarily a density of $P$. In this case  neither $L(\btheta)$ nor $\ell_n(\btheta)$ can be expected to have a unique maximum, not even up to label switching.
Based on a technique inspired by \cite{redner81},  we show that the sequence of constrained ML estimators is asymptotically close to one of the parameters $\btheta^{\star}$ giving the best approximation of $P$ in terms of the KLR.
This technique provides  consistency on the quotient space topology identifying all population log-likelihood maxima.
Define the sets
\[
S(\dot \btheta):=%
\set{%
\btheta \in \Theta_K(P) : %
\int \log \psi(\bx;\btheta)dP(\bx) = \int \log \psi(\bx;\dot \btheta)dP(\bx) %
},
\]
\[
\mathcal{T}(\dot \btheta, \eps):= %
\set{%
\btheta \in \Theta_K(P) : %
\|\btheta - \ddot{\btheta} \| <  \eps \;  \forall \; \ddot \theta \in S(\dot \btheta)%
},
\quad \text{for any}  \quad \eps >0.
\]
Note that $S(\dot \btheta)$ and $\mathcal{T}(\dot \btheta, \eps)$ do not depend on the specific likelihood maximiser if $\dot\btheta=\btheta^\star$.
\begin{theorem}[Consistency]\label{thm:mle_consistency}
Under \Crefrange{asm:no-mass}{asm:monotone-LK},  for every $\eps >0$ and every sequence of maximizers $\btheta_n$ of $\ell_n(\cdot)$:~
$\lim_{n\to\infty} \prob{ \btheta_n \in \mathcal{T}(\btheta^\star, \eps) }  = 1$.
\end{theorem}

\begin{proof} 
Because of Lemmas \ref{lem:bounded-eig}-\ref{lem:bounded-eig-ln} and  \ref{lem:bounded-mu}-\ref{lem:bounded-mu-ln}  there is a compact set $T_{\tilde{\Theta}_K} \subset \tilde{\Theta}_K$ so that all $\btheta_n \in T_{\tilde{\Theta}_K}$ for large enough $n$ a.s.
Using Lemma \ref{lem:bounded-eig}, $|\log \psi(\bx, \btheta) | \leq C$ for some finite constant $C$  for all $\btheta \in T_{\tilde{\Theta}_K}$.
Sufficient conditions for Theorem 2 in \cite{Jennrich69} are satisfied, and therefore  $\sup_{\btheta \in T_{\tilde{\Theta}_K}} \abs{\ell_{n}(\btheta) - L(\btheta)} \conv 0$ $P$-a.s.
Applying the  same argument as in the proof of Theorem 5.7 in \cite{vandervaart00-book},
we obtain   $L({\btheta}_n)  \conv  L({\btheta}^\star)$ $P$-a.s.
By continuity of $L(\cdot)$ and \Cref{thm:asy_existence}  we have that for every $\eps>0$ there exists a $\delta>0$ such that $L(\btheta^\star)-L(\btheta)>\delta$ for all  $\btheta \in T_{\tilde{\Theta}_K}\setminus\mathcal{T}(\btheta^\star, \eps)$.
Denote $(\Omega, \mathcal{A}, P)$ the probability space where the sample random variables are  defined, and consider the following events
\begin{displaymath}
A_n:=\set{\omega \in \Omega : \; \btheta_n  \in T_{\tilde{\Theta}_K}\setminus\mathcal{T}(\btheta^\star,\eps)}, 
\end{displaymath} 
and 
\begin{displaymath}
B_n := \set{\omega \in \Omega : \; L(\btheta^\star) -  L(\btheta_n)> \delta}.
\end{displaymath} 
Clearly $A_n \subseteq B_n$ for all $n$. $P(B_n) \conv 0$ for $n \conv \infty$ implies  $P(A_n) \conv 0$. The latter proves the result.
\end{proof}


\section{The MLE functional in the mixture case}
\label{sec:consistency-2}

After having showed the nonparametric consistency of the MLE for its own
canonical functional, here we investigate the canonical functional for a $P$
that is a mixture of $K$ well enough separated
nonparametric components $Q_1,\ldots,Q_K$ \REV{(see \cite{hunter24} for methods that address this problem explicitly)}. 

Such a $P$ can be interpreted as generating $K$ well separated clusters,
even though the $Q_k,\ k=1,\ldots,K,$ can be chosen so flexibly that one or more
of them themselves could be multimodal or even a mixture of well separated
components generating homogeneous data. As the MLE functional of the mixture of
ESDs enforces $K$ components to be fitted to $P$, it seems reasonable,
interpreting these as corresponding to
$K$ clusters, to expect that they align with
$Q_1,\ldots,Q_K$ also in the latter case if the separation between
$Q_1,\ldots,Q_K$ is stronger than the separation between any ``subcomponents''
of any $Q_k$.

For $P$ of this type, the clusters generated
by the $K$ ESD mixture components of the MLE functional (ESDC)
will indeed approximately
correspond to the ``central regions'' of $Q_1,\ldots,Q_K$, and the parameters
of the ESDC will approximate the parameters of the MLE canonical functionals of
separate ESDs evaluated at each of $Q_1,\ldots,Q_K$ alone. For example, if
the ESD family is chosen as Gaussian with flexible means and covariance
matrices, the corresponding parameters of the ESDC will approximate the
mean and covariance matrix functionals of $Q_1,\ldots,Q_K$.


Here are some definitions and assumptions.
Let $Q_1,\ldots,Q_K$ be distributions on $\R^p$ (generally the same notation refers to distributions and their cumulative distribution functions)
parameterized in such a way that 0 is their ``center'' in some sense; it
could be the mode, the mean, \REV{a (however defined)} multivariate median or quantile; important
is only that $Q_k$ is defined relative to 0.
Let $\xi_1,\ldots,\xi_K>0$ mixture proportions with $\sum_{k=1}^K\xi_k=1$.
For $m\in\N,\ k\in\{1,\ldots,K\}$ let
$\brho_{mk}\in \R^p$ sequences so that
\[
\lim_{m\to\infty}\min_{k_1\neq k_2\in\{1,\ldots,K\}}\|\brho_{mk_1}-\brho_{mk_2}\|
=\infty.
\]
Define a sequence of distributions $P_m$ on $\R^p$ by
$P_m(\bx):=\sum_{k=1}^K\xi_kQ_k(\bx-\brho_{mk})$. 
The mixture $P_m$ is constructed in such a way that
its mixture components for increasing $m$ become better and better separated,
although they are nonparametric and may have non-vanishung densities, so there
may be overlap betwen them even for arbitrarily large $m$.

Consider, for $\epsilon>0$,
the ``central set'' $\{\bx:\ \|\bx\|<\epsilon\}$ about 0. 
$\epsilon$ can be chosen large enough that for arbitrarily small 
$\eta>0$:
\begin{equation}
  \label{eq:epslarge}
\forall k\in\{1,\ldots,K\}:\
Q_k\{\|\bx\|<\epsilon\} \ge 1-\eta,  
\end{equation}
which in particular implies that
\begin{equation}\label{eq:a5}
\exists \delta>0:\
\forall k\in\{1,\ldots,K\}:\
\xi_k Q_k\{\|\bx\|<\epsilon\} \ge \delta.
\end{equation}
The following theorem states that in this setup, when evaluating 
$\btheta^\star(P_m)$, eventually the different clusters 
include
the full (arbitrarily large) central sets of the different mixture components, 
and in this sense
the clustering corresponds to the mixture structure. We require:
\begin{assumption}\label{asm:lambdazero} For any sequence $(\bmu_n,\bSigma_n)_{n\in \N}$, $Q_k,\ k=1,\ldots,K$ for which $\frac{\lmax^*(\bSigma_n)}{\lmin^*(\bSigma_n)}\leq \gamma$:
\begin{displaymath} 
\lim_{n\to\infty} \lmin^*(\bSigma_n)=0\Rightarrow \lim_{n\to\infty} \int \log f(\bx;\bmu_n,\bSigma_n) dQ_k(\bx)=-\infty.   
\end{displaymath}
\end{assumption}
\begin{assumption} \label{asm:intlogg}
$\exists c_0>-\infty$ so that for all $k\in\{1,\ldots,K\}\
:\ \int \log g(\|\bx\|) dQ_k(\bx)\ge c_0.$
\end{assumption}
\begin{remark}\label{rmk:consass}
Assumption \ref{asm:lambdazero} is e.g. fulfilled for $Q_k$ if $Q_k$ is not concentrated on a single point and $g(x)\in o(x^{-\delta})$, where $\delta>\frac{p}{2\epsilon}$, $\epsilon>0$ so that there are two disjoint sets $A$ and $B:\ Q_k(A)\ge\epsilon,\ Q_k(B)\ge\epsilon, \inf_{x\in A, y\in B}\frac{\|x-y\|}{2}=:\eta>0$. This is because in that case, using \cref{eq:esd_eigen},
\begin{eqnarray*}
&\int \log f(\bx;\bmu_n,\bSigma_n) dQ_k(\bx) \le &\\  
&\epsilon\left[\log g\left(\lmax^*(\bSigma_n)^{-1}\eta^2\right)-\frac{p}{2}\log  \lmin^*(\bSigma_n)\right]
+(1-\epsilon)\left[\log g(0)-\frac{p}{2}\log  \lmin^*(\bSigma_n)\right]&\\
&\in o\left((\epsilon\delta-\frac{p}{2}) \log \lmax^*(\bSigma_n)\right).&
\end{eqnarray*}
With similar reasoning, $Q_k$ continuous and  $g(x)\in o(x^{-\delta})$ with $\delta>\frac{p}{2}$ will fulfill Assumption \ref{asm:lambdazero}.

Assumption 
\ref{asm:intlogg} may be violated if $Q_k$ has far heavier tails than $f$; e.g., if $f$ is Gaussian, it amounts to $Q_k$ having an existing covariance matrix.
\end{remark}

\begin{theorem}\label{thm:level} With the above definitions, under Assumption \ref{asm:lambdazero},
for large enough $m$, the clusters of $\btheta^\star(P_m)$ can be numbered 
in such a
way that for $k\in\{1,\ldots,K\}:$
\[
B_\epsilon(\brho_{mk}):=\{\bx:\ \|\bx-\brho_{mk}\|<\epsilon\} \subseteq C_{mk}=\{\bx:\ \cl(\bx,\btheta^\star(P_m))=k\}.
\]
\end{theorem}

\begin{proof} Show the following statements:
\begin{description}
\item[S1] For $\tilde{\btheta}_m$ defined by $\tilde{\pi}_{mk}=\xi_{k},\ 
\tilde{\bmu}_{mk}=\brho_{mk},\ \tilde{\bSigma}_{mk}=\bI_p,\ k=1,\ldots,K$:
\begin{equation}
  \label{eq:s1}
\exists m^->-\infty\ \forall m:\ L(\tilde{\btheta}_m, P_m)\ge m^-.
\end{equation} 
\item[S2] For a sequence $(\btheta_m)_{m\in\N}$ let 
\[
D_m(\btheta_m):=\max_{1\le k\le K}\min_{1\le j\le K}\|\brho_{mk}-\bmu_{mj}\|. 
\]
If $\limsup_{m\to\infty}D_m(\btheta_m)=\infty$, then 
$\liminf_{m\to\infty}L(\btheta_m, P_m)= -\infty$.
\item[S3] The following hold for $\btheta_m:=\btheta^\star(P_m)$: 
There are constants $0<c_1<c_2<\infty$ 
independent of $m$ so that $c_1<\lmin(\btheta_m)<c_2$, and there is a constant $c_3>0$ so that for large enough $m,\ k=1,\ldots,K:\ \pi_{mk}\ge c_3$.
\item[S4]
If $\exists m_D<\infty$ so that $\forall m:\ D_m(\btheta_m)\le m_D$, and S3 holds for $(\btheta_m)_{m\in\N}$,
then for large enough $m$
the components of $\btheta_m$ can be numbered so that for $k=1,\ldots,K$: 
\begin{equation}\label{eq:dm}
B_\epsilon(\brho_{mk}) \subseteq C_{mk}.
\end{equation}
\end{description}
S1 together with S2 imply that 
$\limsup_{m\to\infty}D_m(\btheta^\star(P_m))<\infty$, so that, together with S3, 
$(\btheta^\star(P_m))_{m\in\N}$ fulfills the condition for S4, from which the theorem follows. 
~\\~\\
{\it Proof of S1:}
\begin{eqnarray*}
L(\tilde{\btheta}_m, P_m) &=& \int\log\left(\sum_{j=1}^K \tilde{\pi}_{mj} f(\bx; \tilde{\bmu}_{mj}, \tilde{\bSigma}_{mj})\right) d\sum_{k=1}^K\xi_kQ_k(\bx-\brho_{mk})\\
&= & \sum_{k=1}^K \xi_k  \int\log\left(\sum_{j=1}^K \xi_j f(\bx; \brho_{mj}, \bI_p)\right) dQ_k(\bx-\brho_{mk})\\
&\ge &  \sum_{k=1}^K \xi_k  \int\log\left(\xi_kf(\bx; \brho_{mk}, \bI_p)\right)dQ_k(\bx-\brho_{mk})\\
&=& \sum_{k=1}^K \xi_k  \int\log\left(\xi_kg(\tr{\bx}\bx))\right)dQ_k(\bx)=m^-,
\end{eqnarray*}
independently of $m$. $m^->-\infty$ follows from Assumption \ref{asm:intlogg}.
~\\~\\
{\it Proof of S2:}~\\
It suffices to consider $\lim_{m\to\infty}D_m(\btheta_m)=\infty$ because in case
this holds for the $\limsup$, there is a subsequence diverging to $\infty$. 
W.l.o.g., number mixture components so that $D_m(\btheta_m)=\|\brho_{mk^*}-\bmu_{mk^*}\|$ for a 
fixed $k^*\in\{1,\ldots,K\}$.

The following is required for showing that $\lim_{m\to\infty}L(\btheta_m, P_m)= -\infty$: Because of the properties of $f$ and $g$, 
for $k=1,\ldots,K$ and any sequence $(\bmu_n,\bSigma_n)_{n\in \N}$ for which $\frac{\lmax^*(\bSigma_n)}{\lmin^*(\bSigma_n)}\leq \gamma$:
\begin{equation}\label{eq:lminzeromininfty}
\lim_{n\to\infty}\lmin^*(\bSigma_n)=\infty \Rightarrow 
\int\log f(\bx;\bmu_n,\bSigma_n)dQ_k(\bx)=-\infty.
\end{equation}
Furthermore, 
\begin{equation}
  \label{eq:rhomuinfty}
  \lim_{n\to\infty}|\bmu_n|=\infty \Rightarrow \int\log f(\bx;\bmu_n,\bSigma_n)dQ_k(\bx)=-\infty,
\end{equation}
regardless of whether $\lmin^*(\bSigma_n)$ is bounded (in which case $\log f(\bx;\bmu_n,\bSigma_n)\conv -\infty$ for $\bx\in S$ with $Q_k(S)$ arbitrarily close to 1), diverges to $\infty$ (in which case \cref{eq:lminzeromininfty} obtains), or converges to zero (Assumption \ref{asm:lambdazero}). 

Observe
\begin{equation}\label{eq:lmixq}
  L(\btheta_m,P_m) = \xi_{k^*}\int\log\left(\sum_{j=1}^K \pi_{mj} f(\bx; \bmu_{mj},\bSigma_{mj})\right) dQ_{k^*}(\bx-\brho_{mk^*})+L^*_m,
\end{equation}
where
\begin{displaymath}
L^*_m:=\sum_{k\neq k^*}  \xi_k\int\log\left(\sum_{j=1}^K \pi_{mj} f(\bx; \bmu_{mj},\bSigma_{mj})\right) dQ_k(\bx-\brho_{mk}).
\end{displaymath}
Because of \cref{eq:rhomuinfty}, the first term of \cref{eq:lmixq} diverges to $-\infty$, and so does $L(\btheta_m,P_m)$ if $L^*_m$
is bounded from above. Assumption \ref{asm:lambdazero} implies that 
$\lmin(\btheta_m)\ge c>0$, therefore 
$f(\bx; \bmu_{mk},\bSigma_{mk})\le c^{-p/2}g(0)<\infty$, and $L^*_m\le c^{-p/2}g(0)$.
This proves S2.
~\\~\\
{\it Proof of S3:}\\
$c_1<\lmin(\btheta^\star(P_m))<c_2$ holds by S1, Assumption \ref{asm:lambdazero}, and \cref{eq:lminzeromininfty}. 

Because of S2, $\exists m_D<\infty$ so that $\forall m:\ D_m(\btheta^\star(P_m))\le m_D$. 
Number the components of $\btheta_m$ so that for $k=1,\ldots,K$:
\begin{equation}\label{eq:match}
\|\brho_{mk}-\bmu_{mk}\|=\min_{1\le j\le K}\|\brho_{mk}-\bmu_{mj}\|.
\end{equation}
This is possible because of $D_m(\btheta_m)\le m_D$. 

Assume w.l.o.g. that $\pi_{m1}\conv 0$. Then, \cref{eq:lmixq} holds with $k^*=1$. $L^*_m$ is once more bounded from above, and 
\begin{equation}\label{eq:q1infty}
  \lim_{m\to\infty}  \int\log\left(\sum_{j=1}^K \pi_{mj} f(\bx; \bmu_{mj},\bSigma_{mj})\right)
dQ_{1}(\bx-\brho_{m1})=-\infty,
\end{equation}
because $\pi_{m1}\conv 0$, $f(\bx; \bmu_{m1},\bSigma_{m1})\le c_1^{-p/2}g(0)<\infty$ as in the proof of S2, and for $j\neq 1:\ f(\bx; \bmu_{mj},\bSigma_{mj})\conv 0$ for $\bx\in S$ with 
$\int_S dQ_1(\bx-\brho_{m1})$ arbitrarily close to 1. But for $\btheta_m=\btheta^\star(P_m)$, 
\cref{eq:q1infty} is in contradiction to S1.
~\\~\\
{\it Proof of S4:}~\\
For large enough $m$ number the components  of $\btheta_m$ according to \cref{eq:match}. 
For $\tilde\bx\in B_\epsilon(\bzero)$
so that $\bx=\tilde\bx+\brho_{mk}\in B_\epsilon(\brho_{mk})$ consider
\begin{displaymath}
\tau_{k}(\bx; \btheta_m) = 
\frac{\pi_{mk} f(\bx; \bmu_{mk}, \bSigma_{mk})}
{\sum_{j=1}^K \pi_{mj} f(\bx; \bmu_{mj}, \bSigma_{mj})}.
\end{displaymath}
Because of S3, $\pi_{mk}\ge c_3$, and $f(\bx; \bmu_{mk}, \bSigma_{mk})\ge c_4>0$, whereas
for $j\neq k:\ f(\bx; \bmu_{mj}, \bSigma_{mj}) \allowbreak \conv 0$ uniformly for 
$\tilde\bx\in B_\epsilon(\bzero)$,
so that $\tau_{k}(\bx; \btheta_m)\conv 1$,
and for large enough $m:\ \forall \bx:\ \cl(\bx,\btheta_m)=k$. 
\end{proof}
\begin{remark} It is not essential that $B_\epsilon(\brho_{mk})$ is defined based on the
Euclidean distance; other $L_p$-distances will work as well as Mahalanobis distances with any
fixed covariance matrix.
\end{remark} 

\begin{remark} For $m$ large, the different components $Q_k,\ k=1,\ldots,K,$ are shifted
by $\brho_{mk}$ and their central sets $B_\epsilon(\brho_{mk})$ are therefore very far away from 
each other. It may therefore not seem surprising that these ultimately belong to different 
clusters. But the statement is not trivial. The $Q_k$ can have densities that are nowhere zero
and even have heavy tails (Assumption \ref{asm:intlogg} will then require $f$ to have heavy tails, too), so that there will always be overlap between the $Q_k$. There are clustering methods
with fixed $K$ that for $n\conv\infty$ will not match the $K$ different clusters to $B_\epsilon(\brho_{m1}),\ldots,B_\epsilon(\brho_{mK})$:
\begin{itemize}
\item If the $Q_k$ (or even at least one of them) 
are chosen so that for $n\conv\infty$, the largest distance $D_{max,n}$ 
of an observation to
its nearest neighbour goes to $\infty$ in probability (which holds 
for distributions with heavy enough tails, see Theorem 5.3 in \cite{JamJan15}),
then for any given but arbitrarily large 
$m$, $n$ can be chosen large enough so that $D_{max,n}$ is arbitrarily much larger than 
$\max_{1\le j,k\le K} \|\brho_{mj}-\brho_{mk}\|$ with arbitrarily large probability. Then,
the Single Linkage dendrogram based on the Euclidean distance will merge different central sets earlier than connecting the observation that is farthest from its nearest neighbour with anything else. 
\item A trimmed clustering method trimming a proportion of $\alpha$ observations can trim a
complete central set if $\alpha\ge\xi_k$ for some $k$. Similarly, RIMLE (\cite{ch17-jmlr}) may classify a complete central set as noise.    
\end{itemize}
Arguably there are situations in which the behaviour of trimmed clustering and RIMLE as 
explained above may be seen as desirable, namely if one of the $Q_k$ looks like
modelling more than one cluster; e.g., it can be bimodal with two clearly separated modes. 
It may then be seen as appropriate if this attracts more than one of the $K$ clusters,
whereas another $Q_k$ with either low probability or low (smoothed) density may be classified
as generating outliers. 
\end{remark}
The following theorem states that for the setup of Theorem \ref{thm:level}, the estimators of $\bmu_{mk}$ and $\bSigma_{mk}$ converge to the estimators for the individual mixture components $Q_{mk}$, so that the increasing separation of mixture components for $m\conv\infty$ implies that ultimately the characteristics of $Q_{mk}$ defined in terms of the ESD densities $f$
can be estimated without influence of the other mixture components. For example, if $f$ defines a Gaussian location-scale family, $\bmu^\star_{mk}$ will converge toward the mean of $Q_{mk}$, and $\bSigma^\star_{mk}$ will converge to the covariance matrix of $Q_{mk}$, see Corollary \ref{cor:gaussest}, which shows that the required Assumption \ref{asm:uniquemax} below holds at least in this case.  

For $k=1,\ldots,K$, let
\[
\tilde\bkappa_k:=(\tilde\bmu_k,\tilde\bSigma_k)=\argmax_{\bkappa}\tilde L(\bkappa, Q_k),\
\tilde L(\bkappa, Q):=
\int \log f(\bx;\bkappa)dQ_k(\bx).
\]
The corresponding functionals for $Q_{mk}=Q_k(\bullet-\brho_{mk})$ are
\begin{equation}
  \label{eq:qest}
  \tilde\bmu_{mk}:=\tilde\bmu_k+\brho_{mk},\ \tilde\bSigma_{mk}:=\tilde\bSigma_{k}.
\end{equation}
\begin{assumption}\label{asm:uniquemax} For given $Q_k$, 
\[ 
\forall\eps>0\ \exists\beta>0:\ \|\bkappa-\tilde\bkappa_k\|>\eps \Rightarrow
L(\tilde\bkappa_k,Q_k)-L(\bkappa,Q_k)>\beta.
\]
\end{assumption}
\begin{assumption}\label{asm:ksigmalambda} For $\tilde\btheta=(\xi_1,\ldots,\xi_k,\tilde\bkappa_1,\ldots,\tilde\bkappa_K)$: $\frac{\lmin(\tilde\btheta)}{\lmax(\tilde\btheta)}\le \gamma$.
\end{assumption}  
 
\begin{theorem}\label{thm:estimatorsconverge}  With the above definitions, under Assumptions \ref{asm:lambdazero} and \ref{asm:intlogg}, 
for large enough $m$, the clusters of $\btheta^\star(P_m)$ can be numbered 
in such a
way that for $k\in\{1,\ldots,K\},$ with $\btheta^\star_m:=\btheta^\star(P_m)$.
\begin{equation}
  \label{eq:piconv}
  \lim_{m\to\infty}\|\pi^\star_{mk}-\xi_{k}\|=0,
\end{equation}
and if $Q_1,\ldots,Q_K$ fulfill Assumption \ref{asm:ksigmalambda}, then for those $Q_k$ that fulfill Assumption \ref{asm:uniquemax}:
\begin{equation}\label{eq:kappaconv}
  \lim_{m\to\infty}\|\bkappa^\star_{mk}-\tilde\bkappa_{mk}\|=0.
\end{equation}
\end{theorem} 
\begin{proof} Assume throughout that mixture components are numbered according to \cref{eq:match}. Show the following statements:
\begin{description}
\item[S1] For $\btheta^\star_m=\btheta^\star(P_m)$, using the notation of 
\cref{eq:tau}, and $\btheta^\star_{mk}$ denoting the parameters in $\btheta^\star_m$ belonging to mixture component $k,\ k=1,\ldots,K$,  
\[
\btheta^\star_{mk}=\argmax_{\btheta}\int \tau_k(\bx;\btheta^\star_{mk})(\log \pi_k+\log 
f(\bx;\bkappa_k)) dP_m(\bx).
\]
\item[S2] \cref{eq:piconv} follows from S1.
\item[S3] From S1,
\[
\bkappa^\star_{mk}=\argmax_{\bkappa}q_{mk}(\bkappa),
\]
where $q_{mk}(\bkappa)$ is a function so that 
$q_{mk}(\bkappa)-\int \log f(\bx;\bkappa)dQ_k(\bx-\brho_{mk})$ converges
uniformly to 0.  
\item[S4] \cref{eq:kappaconv} follows from S3 and Assumption \ref{asm:uniquemax}.
\end{description}
{\it Proof of S1:}~\\ The proof is based on the two-step form of a mixture model involving for each observed $X$ an unobserved random variable $\zeta$ indicating the mixture component that has generated $X$, see Section \ref{sec:mix-clustering-classification}. S1 is the population version of the fixed point equation on which the EM-algorithm is based, see Section 3 and in particular Corollary 2 of \cite{delaru77}. 

Regarding the mixture \cref{eq:psi}, on which the ML-estimation is based, let
$\tilde \psi(\bullet;\btheta)$ be the joint density of $Y=(X,\zeta)$, and
$\bar\psi(\bullet|\bx;\btheta)$ be the conditional density of $Y$ given $X=\bx$
so that for all $\by,\ \bx$:
\begin{equation}
  \label{eq:krat}
  \bar\psi(\by|\bx;\btheta)=\frac{\tilde \psi(\by;\btheta)}{\psi(\bx;\btheta)}.
\end{equation}
Define
\begin{eqnarray*}
G(\btheta'|\btheta)&:=&\int\int \log(\tilde \psi(\by;\btheta'))\bar\psi(\by|\bx;\btheta)d\by dP(\bx)\\
H(\btheta'|\btheta)&:=&\int\int\log(\bar\psi(\by|\bx;\btheta'))\bar\psi(\by|\bx;\btheta)d\by dP(\bx).
\end{eqnarray*}
Then, following \cite{delaru77},  
\cref{eq:krat} still holds averaged over the $y|\bx$ assumed distributed 
according to $\bar\psi(\by|\bx;\btheta)$, and then averaging over $P_m(\bx)$ yields
\begin{eqnarray*}
 L(\btheta',P)= G(\btheta'|\btheta)-H(\btheta'|\btheta).
\end{eqnarray*}
For given $\bx$, formula (1e6.6) in \cite{rao65} implies
\begin{eqnarray*}
  \int (\log \bar\psi(\by|\bx;\btheta)-\log \bar\psi(\by|\bx;\btheta'))\bar\psi(\by|\bx;\btheta)dy&\ge& 0\\
\Rightarrow H(\btheta|\btheta)&\ge& H(\btheta'|\btheta).
\end{eqnarray*}
Let $M(\btheta^\star_m):=\argmax_{\btheta} G(\btheta|\btheta^\star_m)$. Then
\[
L(M(\btheta^\star_m))=Q(M(\btheta^\star_m)|\btheta^\star_m)-H(M(\btheta^\star_m)|\btheta^\star_m)\ge G(\btheta^\star_m|\btheta^\star_m)-H(\btheta^\star_m|\btheta^\star_m)=L(\btheta^\star_m),
\]
implying $\btheta^\star_m=\argmax_{\btheta} G(\btheta|\btheta^\star_m)$ (otherwise the ``$\ge$'' above would be ``$>$'' and $M(\btheta^\star_m)$ would improve $L$). Note that
\begin{displaymath}
  G(\btheta|\btheta^\star_m)=\int \sum_{k=1}^K \tau_k(\bx;\btheta^\star_{mk})(\log \pi_k+\log 
f(\bx;\bkappa_k)) dP_m(\bx),
\end{displaymath}
which can be maximised for each $k$ separately, leading to S1.
~\\~\\
{\it Proof of S2:}~\\  
$\int \tau_k(\bx;\btheta^\star_{mk})(\log \pi_k+\log 
f(\bx;\bkappa_k)) dP_m(\bx)$ can be maximised separately over $\pi_1,\ldots,\pi_K$ and $\bkappa_1,\ldots,\bkappa_K$. Maximising $\int \tau_k(\bx;\btheta^\star_{mk})\log \pi_kdP_m(\bx)$ yields 
\begin{equation}\label{eq:bepsdecomp}
\pi^\star_{mk}=\int\tau_k(\bx;\btheta^\star_{mk})dP_m(\bx)=
\int_{B_\eps(\brho_{mk})}\tau_k(\bx;\btheta^\star_{mk})dP_m(\bx)+
\int_{B_\eps(\brho_{mk})^c}\tau_k(\bx;\btheta^\star_{mk})dP_m(\bx),
\end{equation}
$B_\eps(\brho_{mk})$ chosen as in the proof of S4 of Theorem \ref{thm:level}.
From there, $\tau_k(\bx;\btheta^\star_{mk})\conv 1$ 
for $\tilde\bx\in B_\epsilon(\bzero)$, i.e., $\bx\in B_\epsilon(\brho_{mk})$,
therefore
\begin{equation}\label{eq:tauto1}
\lim_{m\to\infty}
  \left|\int_{B_\eps(\brho_{mk})}\tau_k(\bx;\btheta^\star_{mk})dP_m(\bx)-
\int_{B_\eps(\brho_{mk})}d\sum_{k=1}^K\xi_kQ_{mk}(\bx)\right|= 0,
\end{equation}
and furthermore, for $j\neq k,\ m\conv\infty$, because of \cref{eq:epslarge},
\begin{equation}\label{eq:outsidevanish}
Q_{mj}(B_\eps(\brho_{mk}))\conv 0 \Rightarrow
  \int_{B_\eps(\brho_{mk})}d\xi_jQ_{mj}(\bx)\conv 0.
\end{equation}
$\eps$ can be chosen so large that $Q_{mk}(B_\eps(\brho_{mk}))$ is arbitrarily 
close to 1, and due to \cref{eq:outsidevanish}, 
$\int_{B_\eps(\brho_{mk})}d\sum_{k=1}^K\xi_kQ_{mk}(\bx)$ is arbitrarily close to
$\xi_k$. Furthermore, assuming $m$ large enough that 
$B_\eps(\brho_{mj}),\ j=1,\ldots,K$, do not intersect, with $\tilde B=\left(\bigcup_{j=1}^K B_\eps(\brho_{mj})\right)^c$,
\begin{displaymath}
 \int_{B_\eps(\brho_{mk})^c}\tau_k(\bx;\btheta^\star_{mk})dP_m(\bx)= 
\sum_{j\neq k}\int_{B_\eps(\brho_{mj})}\tau_k(\bx;\btheta^\star_{mk})dP_m(\bx)+
\int_{\tilde B}\tau_k(\bx;\btheta^\star_{mk})dP_m(\bx).
\end{displaymath}
For the same reason as \cref{eq:tauto1} and \cref{eq:outsidevanish} (with 
roles of $j$ and $k$ inverted),
\begin{displaymath}
 \sum_{j\neq k}\int_{B_\eps(\brho_{mj})}\tau_k(\bx;\btheta^\star_{mk})dP_m(\bx)\conv 0. 
\end{displaymath}
Because of \cref{eq:epslarge}, $P_m(\tilde B)$ is arbitrarily small for $\eps$
large enough, as is
$\int_{\tilde B}\tau_k(\bx;\btheta^\star_{mk})dP_m(\bx)$. Putting everything 
together, \cref{eq:bepsdecomp} implies \cref{eq:piconv}.
~\\~\\
{\it Proof of S3:}~\\ From S1, $\bkappa^\star_{mk}=\argmax_{\bkappa}\tilde q_{mk}(\bkappa)$
with
\[
\tilde q_{mk}(\bkappa)=\int \tau_k(\bx;\btheta^\star_{mk})\log 
f(\bx;\bkappa_k) dP_m(\bx)=\int \tau_k(\bx;\btheta^\star_{mk})\log 
f(\bx;\bkappa_k) d\sum_{k=1}^K\xi_kQ_{mk}(\bx).
\]
As in \cref{eq:bepsdecomp} and the proof of S2,
\begin{displaymath}
 \tilde q_{mk}(\bkappa)=\int_{B_\eps(\brho_{mk})} \tau_k(\bx;\btheta^\star_{mk})\log 
f(\bx;\bkappa_k) dP_m(\bx) +\int_{B_\eps(\brho_{mk})^c} \tau_k(\bx;\btheta^\star_{mk})\log 
f(\bx;\bkappa_k) dP_m(\bx),
\end{displaymath}
and
\begin{eqnarray*}
\int_{B_\eps(\brho_{mk})^c} \tau_k(\bx;\btheta^\star_{mk})\log 
f(\bx;\bkappa_k) dP_m(\bx)&=&\sum_{j\neq k}\int_{B_\eps(\brho_{mj})} \tau_k(\bx;\btheta^\star_{mk})\log 
f(\bx;\bkappa_k) dP_m(\bx)\\
&+&\int_{\tilde B} \tau_k(\bx;\btheta^\star_{mk})\log 
f(\bx;\bkappa_k) dP_m(\bx).
\end{eqnarray*}
Because of Assumption \ref{asm:lambdazero}, regarding $\argmax_{\bkappa}q_{mk}(\bkappa)$, it suffices to consider $\bkappa$ with $\lmin^*(\bSigma)>c_1$ for suitable $c_1>0,\ c_2=\log (c_1^{-p/2}g(0))$, so that $\log 
f(\bx;\bkappa_k)< c_2$. On $B_\eps(\brho_{mk})$, according to the proof of Theorem \ref{thm:level}, $\tau_k(\bx;\btheta^\star_{mk})\conv 1$ and, for $j\neq k$, 
$\tau_j(\bx;\btheta^\star_{mj})\conv 0$, and furthermore, for $m\conv\infty$, 
\begin{displaymath}
  \int_{B_\eps(\brho_{mk})} \tau_k(\bx;\btheta^\star_{mk})\log 
f(\bx;\bkappa_k) dQ_{mj}(\bx) < c_2\int_{B_\eps(\brho_{mk})} \tau_k(\bx;\btheta^\star_{mk})dQ_{mj}(\bx)\conv 0.
\end{displaymath}
Therefore,
\begin{equation}\label{eq:bepslim}
  \lim_{m\to\infty}\left|\int_{B_\eps(\brho_{mk})} \tau_k(\bx;\btheta^\star_{mk})\log 
f(\bx;\bkappa_k) dP_{m}(\bx) -\int_{B_\eps(\brho_{mk})} \log 
f(\bx;\bkappa_k) d\xi_kQ_{mk}(\bx)\right|=0.
\end{equation}
For $j\neq k$, on $B_\eps(\brho_{mj}):\ \tau_k(\bx;\btheta^\star_{mk})\conv 0$,
therefore
\begin{equation}\label{eq:bepslim2}
 \int_{B_\eps(\brho_{mj})} \tau_k(\bx;\btheta^\star_{mk})\log 
f(\bx;\bkappa_k) dP_{m}(\bx)< c_2 \int_{B_\eps(\brho_{mj})} \tau_k(\bx;\btheta^\star_{mk}) dP_{m}(\bx)\conv 0,
\end{equation}
and
\begin{equation}\label{eq:bepslim3}
  \int_{\tilde B} \tau_k(\bx;\btheta^\star_{mk})\log 
f(\bx;\bkappa_k) dP_{m}(\bx)< c_2 P_m(\tilde B).
\end{equation}
Consider $q_{mk}(\bkappa)=\frac{\tilde q_{mk}(\bkappa)}{\xi_k}$ in order to drop 
$\xi_k$ from \cref{eq:bepslim}.
Once more $\eps$ can be chosen large enough that $P_m(\tilde B)$ becomes arbitrarily small, and $Q_{mk}(B_\eps(\brho_{mk}))$ becomes arbitrarily large, so that
\cref{eq:bepslim}, \cref{eq:bepslim2}, and \cref{eq:bepslim3} together imply 
that 
\[
q_{mk}(\bkappa)-\int \log f(\bx;\bkappa)dQ_{mk}(\bx)\conv 0,
\]
and the convergence is uniform over $\bkappa$ as neither 
$\tau_j(\bx;\btheta^\star_{mj})$ for any $j$ nor $\epsilon$ nor $c_2$ 
depend on $\bkappa$.
~\\~\\
{\it Proof of S4:}~\\ Given arbitrarily small $\eps$ and the corresponding $\beta$ from Assumption \ref{asm:uniquemax}, according to S3, 
for $m$ large enough, uniformly
\[
|q_{mk}(\bkappa)-\tilde L(\bkappa,Q_k)|<\frac{\beta}{2}.
\]
This means that $\bkappa$ with $\|\bkappa-\tilde\bkappa_{mk}\|>\epsilon$ cannot
maximise $q_{mk}(\bkappa)$, and therefore $\bkappa^\star_{mk}=\argmax_{\bkappa}q_{mk}(\bkappa)$ will
become arbitrarily close to $\tilde\bkappa_{mk}$, because Assumption \ref{asm:ksigmalambda} enforces $\frac{\lmin(\tilde\btheta)}{\lmax(\tilde\btheta)}<\gamma$, which is also required to hold for $\btheta^\star(P^m)$. This proves \cref{eq:kappaconv}.
\end{proof}

\begin{corollary}\label{cor:gaussest} In the situation of Theorem \ref{thm:estimatorsconverge}, requiring Assumptions \ref{asm:lambdazero} and \ref{asm:intlogg},
if $g$ is chosen so that $f(\bullet;\bmu,\bSigma)$ is the density of a $p$-variate Gaussian distribution with mean $\bmu$ and covariance matrix $\bSigma$, then
\begin{equation}\label{eq:normalmean}
\lim_{m\to\infty}\left\|\bmu^\star_{mk}-\int \bx dQ_k(\bx)-\brho_{mk}\right\|=0.
\end{equation}
If additionally Assumption \ref{asm:ksigmalambda} holds, then
\begin{displaymath}
\lim_{m\to\infty}\left\|\bSigma^\star_{mk}-\int (\bx-\tilde\bmu_k)\tr{(\bx-\tilde\bmu_k)} dQ_k(\bx)\right\|=0.
\end{displaymath}
\end{corollary}
In this case, Assumption \ref{asm:lambdazero} will hold if $Q_k$ is not concentrated on a single point, see Remark \ref{rmk:consass}. Assumption \ref{asm:intlogg} amounts to $\int\|\tr{\bx}\bx\|dQ_k(\bx)<\infty$.\\

\begin{proof}
Due to Assumption \ref{asm:lambdazero}, consider 
$\bkappa$ with $\lmin^*(\bSigma_k)>c_1$ for $k=1,\ldots,K,\ c_1>0$.
For the multivariate Gaussian (see, e.g., \cite{andolk85})
\begin{eqnarray*}
   \log f(\bx;\bmu,\bSigma)&=& c-\frac{1}{2}\log|\bSigma|-\frac{1}{2}\tr{(\bx-\bmu)}\inv{\bSigma}(\bx-\bmu),\\
\frac{\partial}{\partial \bmu}\log f(\bx;\bmu,\bSigma)&=& \inv{\bSigma}(\bmu-\bx),\\
\frac{\partial}{\partial \inv{\bSigma}}\log f(\bx;\bmu,\bSigma)&=& \frac{1}{2}
\bSigma-\frac{1}{2}(\bx-\mu)\tr{(\bx-\mu)}.
\end{eqnarray*}
For $\bmu \in B_\eps(\bmu_0)$ for any $\bmu_0$, 
\[
\left\|\frac{\partial}{\partial \bmu}\log f(\bx;\bmu,\bSigma)\right\|\le
c_1^{-p}(\|\bmu_0-\bx\|+\eps).
\]
For $\inv{\bSigma}\in B_\eps(\inv{\bSigma_0})$ for any $\bSigma_0$,
\[
\left\|\frac{\partial}{\partial \inv{\bSigma}}\log f(\bx;\bmu,\bSigma)\right\|
\le \frac{1}{2}\left(\bSigma_0+\eps-(\bx-\mu)\tr{(\bx-\mu)}\right).
\]
Both these bounding functions are integrable w.r.t. $P_m$ because of 
Assumption \ref{asm:intlogg}. This means that $\log f(\bx;\bmu,\bSigma)$ can
be differentiated under the integral sign. Standard algebra yields
\[
\tilde\bmu_{mk}=\int x dQ_{mk}(x),\ 
\tilde\bSigma_{mk}=\int (\bx-\tilde\bmu_{mk})\tr{(\bx-\tilde\bmu_{mk})} dQ_{mk}(\bx).
\] 
Assumption \ref{asm:uniquemax} follows because $\log f(\bx;\bmu,\bSigma)$ is 
strictly concave in $\bmu$ and $\inv{\bSigma}$ (\cite{andolk85}), and the
Corollary follows from Theorem \ref{thm:estimatorsconverge}. Assumption \ref{asm:ksigmalambda} is not required for \cref{eq:normalmean}, because $\int x dQ_{mk}(x)$ maximises the Gaussian likelihood regardless of $\bSigma$, cp. S3 and S4 in the proof of Theorem \ref{thm:estimatorsconverge}.
\end{proof} 


\section{Numerical experiments}
\label{sec:simulation}

\begin{figure}[!t]
\centering
\includegraphics[width=\textwidth,height=0.3\textheight,keepaspectratio]{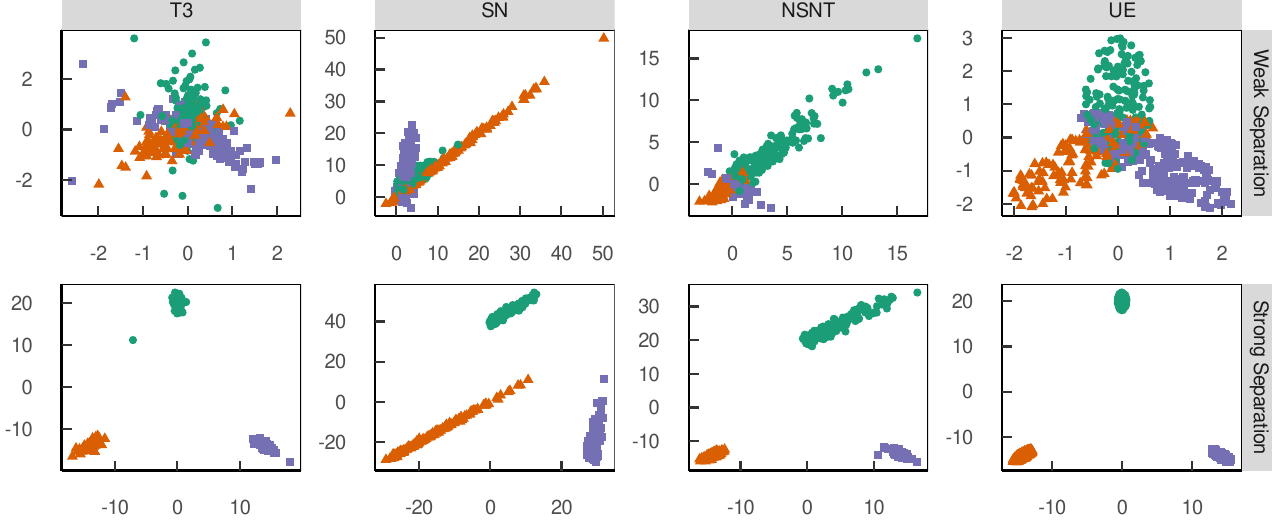}
\caption{Data sets generated from mixture models \textsc{t3, sn}, and \textsc{nsnt}.
``\emph{Weak Separation}'' is obtained setting $s=0.5$ for \textsc{t3} and \textsc{nsnt}, and $s=2$ for \textsc{sn}.
``\emph{Strong Separation}'' is obtained setting $s=20$ for \textsc{t3} and \textsc{nsnt}, and $s=40$ for \textsc{sn}.
The sample size is set to $n=1000$ in all the six cases.}
\label{fig:mc-datasets}
\end{figure}

In order to illustrate the results in Section \ref{sec:consistency-2}, we present a small simulation study. We computed MLEs for mixtures of two  different families of ESDs, namely Gaussian (\textsc{mle-n}) and multivariate $t_5$ (\textsc{mle-t5}), and applied them to mixtures with components that deviate from those assumed.
We generated data from four different mixture models, each with $K=3$, namely a mixture of three bivariate $t_3$-distributions, a mixture of three skew normal distributions (\cite{LeMcL13}), a mixture of one skew normal, one Gaussian and one $t_3$-distribution, \REV{and a mixture of uniform distributions on ellipsoids.}
The previous models are labeled \textsc{t3}, \textsc{sn},  \textsc{nsnt}, and \REV{\textsc{ue}} respectively.
For each model we varied the amount of separation between mixture components. Parameters were chosen as follows. For all models $\pi_1= \pi_2= \pi_3 = 1/3$. Given a separation parameter $s>0$, let
\[
\bmu_1=\tr{(0,s)},\quad 
\bmu_2=\tr{(-s/\sqrt{2},-s/\sqrt{2})},\quad
\bmu_3=\tr{(s/\sqrt{2},-s/\sqrt{2})}, 
\]
\[
\bSigma_1=\left(\begin{array}{cc} 0.1 & 0 \\ 0 & 1 \end{array}\right),\quad 
\bSigma_2=\left(\begin{array}{cc} 0.55 & 0.45 \\ 0.45 & 0.55 \end{array}\right),\quad 
\bSigma_3=\left(\begin{array}{cc} ~~0.55 & -0.45 \\ -0.45 &  ~~0.55 \end{array}\right),
\]
\[
\bdelta_1=\tr{(5,5)},\ \bdelta_2=\tr{(15,15)},\ \bdelta_3=\tr{(1,10)}.
\]
In the \textsc{t3} model, the three mixture components had means $\bmu_1,\, \bmu_2,\,
\bmu_3$ and covariance matrices $\bSigma_1,\, \bSigma_2,\, \bSigma_3$. $s$ was chosen to range from 0.5 to 20.
Skew normal distributions in the \textsc{sn} model were generated according to $\bmu_k+\bdelta_k|Z_0|+Z_1$, where $Z_0\sim{\cal N}_1(0, 1),\ Z_1\sim{\cal N}_2(\bzero,\bSigma_k),\ k=1,2,3$, where ${\cal N}_p$ denotes the $p$-dimensional Gaussian distribution.
If the $k$-th component of the mixture has a skew normal distribution, its expectation is $\bmu_k + \bdelta_k \sqrt{2/\pi}$, and its covariance matrix is $\bSigma_k + \bdelta_k\tr{\bdelta_k} (1-2/\pi)$.
For this model, the separation parameter $s$ was chosen from the range $s=2$ to $s=40$.
For the mixed mixture model \textsc{nsnt}, the skew normal component was generated using the parameters $\bmu_1,\ \bSigma_1,\ \bdelta_1$, the Gaussian component was generated using the mean $\bmu_2$ and the covariance matrix $\bSigma_2$, and the $t_3$-component was generated using the mean $\bmu_3$ and the covariance matrix $\bSigma_3$. For this model, $s$ was chosen in the interval [1,20].
\REV{%
The sampling for the \textsc{ue} model is as follows.
Let $U$ be a uniform distribution generating points on a circle centered at the origin $\bzero$ with radius $r$. $\E{U}=\bzero$, and the radius is such that $\Var{U} = \bI_2$. The $k$th component is defined by the affine transformation $U_k = \bOmega_k U + \bmu_k$, where $\bOmega_k$ is the lower triangular matrix obtained from the Cholesky decomposition of $\bSigma_k$. 
Note that $\E{U_k} = \bmu_k$, and $\Var{U_k}=\bSigma_k$. 
This sampling produces uniformly distributed points on the ellipsoids determined by the previously defined mean and covariance parameters.
As for the previous sampling design, we also consider $s \in [0.5, 20]$.
For each of the four models, we consider a grid of five (logarithmically equispaced) values of $s$ for a total of 15 data-generating processes.
We consider sample sizes $n \in \{250, 1000, 10000\}$ for a total of 60 sample designs.
}
Figure \ref{fig:mc-datasets} shows examples of datasets with $n=1000$ generated by the four models, with two different choices of the separation parameter $s$.

\begin{figure}[!t]
\centering
\includegraphics[width=\textwidth,height=0.33\textheight,keepaspectratio]{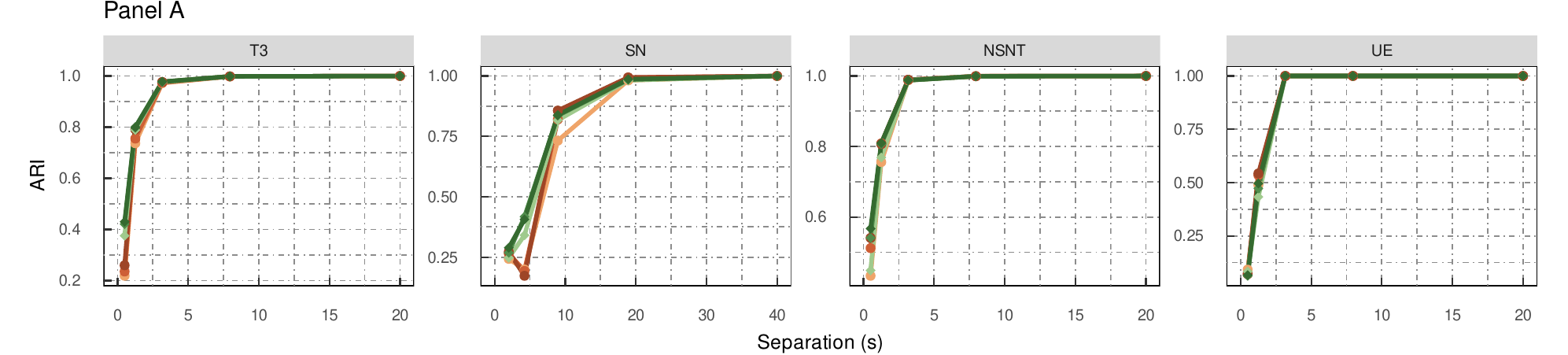}
\includegraphics[width=\textwidth,height=0.33\textheight,keepaspectratio]{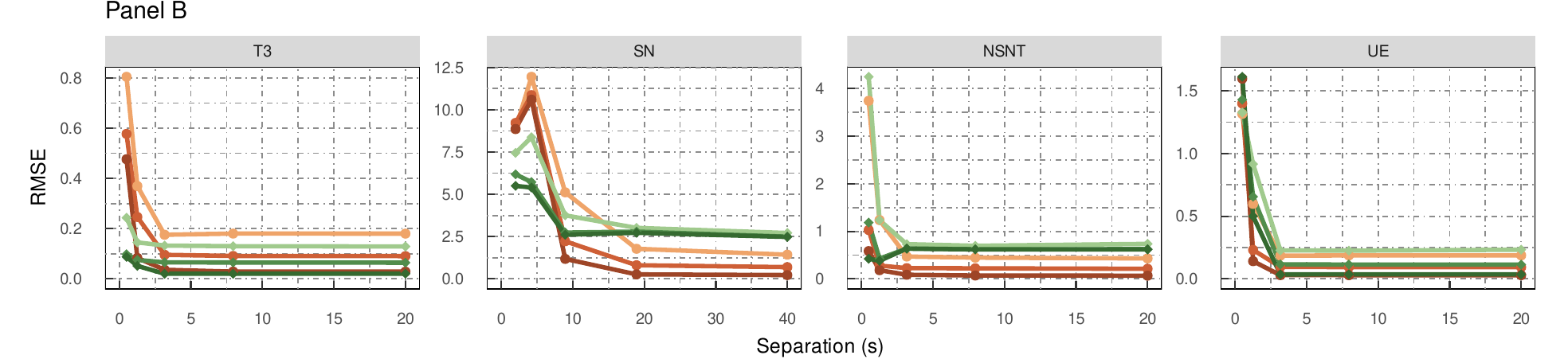}
\includegraphics[width=\textwidth,height=0.33\textheight,keepaspectratio]{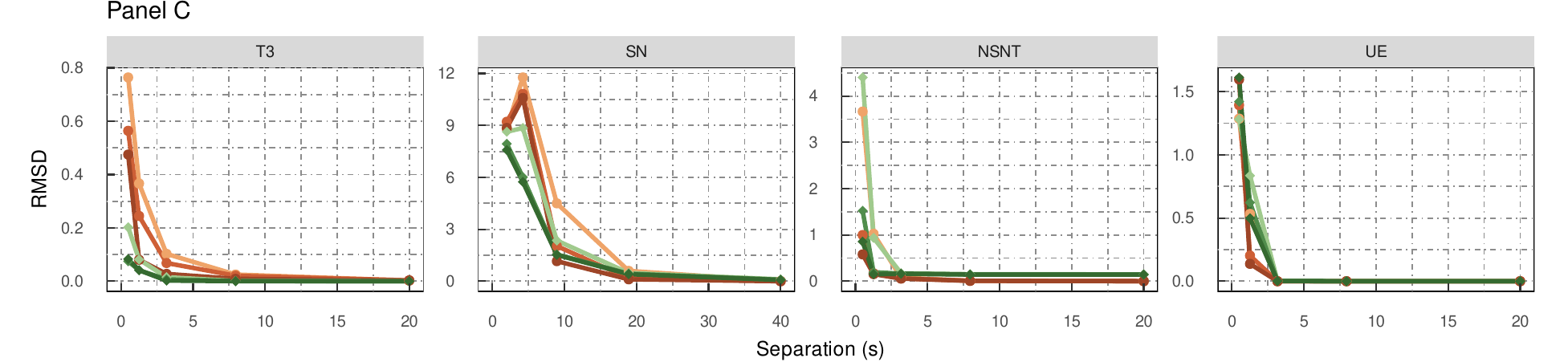}
\includegraphics[width=\textwidth,height=0.1\textheight,keepaspectratio]{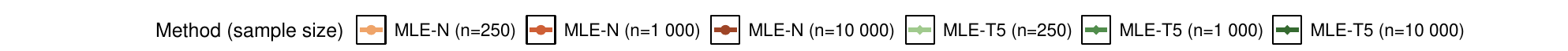}
\caption{Monte Carlo averages for ARI (Panel A), RMSE for mean parameters (Panel B) and RMSD for ML mean functionals (Panel C).}
\label{fig:mc-results}
\end{figure}

We ran $1000$ Monte Carlo replicates with i.i.d. sampling.
For each sample, we computed \textsc{mle-n} and \textsc{mle-t5} with $\gamma = 100$.
After computing the MLEs, the observations were assigned to clusters using the MAP rule.
These were compared with the true mixture components using the Adjusted Rand Index (ARI; \cite{ha85}).
Results are shown in \figref{fig:mc-results}. 
Panel~A shows the ARI averages over the Monte Carlo runs.
Panel~B shows the root mean square error (RMSE) of the mean parameters fitted by \textsc{mln-n} and \textsc{mle-t5} compared to the mean of the generating mixture components. 
Panel~C shows the root mean square deviation (RMSD) between the mean parameters estimated by the fit of \textsc{mle-n} and \textsc{mle-t5}, respectively, and the corresponding estimators computed on only the observations from the true mixture components, i.e., the sample mean (as to be compared to the \textsc{mle-n} component mean estimators), and the ML estimator of location of a $t_5$-distribution (as to be compared to the \textsc{mle-t5} component location estimators). The latter would indicate the best performance one could expect, using true component information, estimating $\tilde\bkappa_{mk}$ from \cref{eq:kappaconv} (Theorem \ref{thm:estimatorsconverge}) by ML.

The results show that with low separation $s$ estimation and clustering do not work well. For the larger values of $s$, the ARI becomes almost perfect, which means that the clustering from the misspecified mixture model matches the true component memberships. This illustrates Theorem \ref{thm:level}. %
\REV{Even with the small sample size $n=250$ already such a good classification is achieved that any improvement with larger $n$ is not visible.
The \textsc{ue} model achieves a performance in line with the other three models, even though here the underlying uniform structure is in some distance from the  Gaussian and Student-t density models on which the estimators \textsc{mle-n} and \textsc{mle-t5} are built.
On the other hand, here there exists $s_0$ such that for $s \geq s_0$, there is a total absence of overlap between the ellipsoidal clustered regions. Note that for the \textsc{ue} model, the tails are as dense as in the center of the distribution. Therefore, we expect a greater overlap for low values of $s$ than the other models considered here. In contrast, the overlap sharply disappears for larger $s \geq s_0$.}
The influence of growing $n$ becomes clear looking at Panel~B, where the 
\REV{improvement of the RMSE for larger values of $n$  is obvious especially in cases of low separation.}
For \textsc{mle-n} the RMSE for $n=10000$ is close to zero for all three setups.
We are interested in particular in examining the RMSE for estimating the mean 
with the aim of illustrating Corollary \ref{cor:gaussest}. For 
$t_3$-distributions, the mean is the center of symmetry, and therefore the
ML-estimator that treats the data as $t_5$ estimates the $t_3$-mean as well, 
and does this better (with lower variance) than the arithmetic mean. 
Correspondingly, \textsc{mle-t5} achieves a lower RMSE than 
\textsc{mle-n} for the model \textsc{t3}, but both will converge to zero for $n\to\infty$. The models \textsc{sn} and \textsc{nsnt} involve asymmetric mixture components for which the ML location functional based on the $t_5$-distribution is not the same as the mean. For this reason, the RMSE for \textsc{mle-t5} comparing it to the mean will not converge to zero for $n\to\infty$, and consequently it seems to hit a nonzero floor in Panel~B.  
\REV{For the uniform model \textsc{ue} the performance is better except for low $s$.}

Panel~C shows that for strong enough separation the location estimators from 
\textsc{mle-n} and \textsc{mle-t5} become indistinguishable from the corresponding estimators computed based on the true component information for models \textsc{t3} and \textsc{sn}, as should be expected from Theorem  \ref{thm:estimatorsconverge} and Panel~A. For model \textsc{nsnt}, however, the RMSD of \textsc{mle-t5} does not seem to converge to zero. Surprised by this, we found numerically that the ML functional based on the $t_5$-distribution computed for the components separately violates  $\frac{\lmin(\tilde\btheta)}{\lmax(\tilde\btheta)}\le \gamma=100$, i.e.,  Assumption \ref{asm:ksigmalambda}, and therefore it cannot be recovered even asymptotically by \textsc{mle-t5}. In fact, also the Gaussian ML functional violates $\frac{\lmin(\tilde\btheta)}{\lmax(\tilde\btheta)}\le 100$ here, but this does not affect the consistency of the means from \textsc{mle-n}, as Assumption \ref{asm:ksigmalambda} is not required for \cref{eq:normalmean}.
\REV{For the uniform model \textsc{ue} the convergence is pretty fast once  $s$ touches the point where the ellipsoids become reasonably separated.}

\section{Conclusion}
\label{sec:conclusions}
Statistical model assumptions are not normally fulfilled in real life situations, and it is of interest how statistical methodology behaves in cases in which the model assumptions are not fulfilled. Consistency for the canonical functional means that for large data sets the estimator will stabilize even if its distributional assumptions are not fulfilled, although still assuming i.i.d. data. Results of this kind can be found in many places. Additionally here we look at specific underlying distributions $P$ that are of such a kind that using the mixture MLE could be of interest for clustering, namely a mixture with $K$ well separated components that can have a different distributional shape from what is assumed. The components of the canonical functional will then correspond in a well defined sense to the components of $P$, although potentially very strong separation is required. Without requiring strong enough separation, such a result cannot be had, and actually it would not be desirable, as the components $Q_k$ of $P$ can themselves be heterogeneous. If different $Q_k$ are not well separated, from the point of view of clustering it may be more sensible to split up a heterogeneous, potentially multimodal, $Q_k$ into subgroups rather than separate it from a $Q_l\neq Q_k$ from which it is not separated. That said, investigating the canonical functional in more general situations is certainly of interest. Furthermore, performance guarantees for fixed $n$ would be welcome, although these will require more restrictive assumptions.

The presented results concern the global optimum of the likelihood function whereas existing algorithms are not guaranteed to find it. Similar results for local optima as found by the EM-algorithm would also be a desirable extension.


\appendix
\section{Affine equivariance of estimators on sphered data}\label{sec:appendix-affine}
\REV{In many situations, an  estimator that is not affine equivariant can be modified into an affine equivariant estimator by being applied on sphered data. We are not aware of any explicit statement of this in the literature (although some statisticians may think that this is folklore knowledge), and we show here a version of this claim that applies in the situation discussed in the present paper.}

\REV{Let $X\in \R^{p\times n}$ a data matrix with $n$ observations and $p$ variables. Let $L:\ \R^{p\times n}\mapsto\R^p$ be an affine equivariant location estimator so that, for $A\in\R^{p\times p}$ invertible and $a\in\R^p$: $L(AX+a)=AL(X)+a$. Let $S:\ \R^{p\times n}\mapsto\R^{p\times p}$ be an affine equivariant scatter estimator so that $S(AX+a)=AS(X)\tr{A}$. Assuming that $S(X)^{-1}$ exists, define $B(X)\in\R^{p\times p}$ by $S(X)^{-1}=\tr{B(X)}B(X)$. Let $P(X)=B(X)(X-L(X))$ be the sphered version of $X$.} 

\REV{Let $T:\ \R^{p\times n}\mapsto \R^K\times \R^{K\times p} \times \R^{K\times p\times p}$ be an estimator of the parameters of an ESD mixture with $K$ components, i.e., $T(X)=(T_1(X),T_2(X),T_3(X))$, where $T_1(X)$ estimates the $K$ proportion parameters $\pi_k$, $T_2(X)$ estimates the $K$ location vectors $\bmu_k$ and, and $T_3(X)$ estimates the $K$ scatter matrices $\bSigma_k,\ k=1,\ldots,K$.}

\REV{$T$ is called affine equivariant if, for $A$ and $a$ as above,
$$
T(AX+a)=(T_1(X),AT_2(X)+a,AT_3(X)\tr{A}).
$$  
We abuse notation here to some extent by not writing out separately what happens to every mixture component, i.e., by $AT_2(X)+a$ we mean that the transformation $L\mapsto AL+a$ is applied to all $K$ component location estimators etc.
\begin{proposition} For any not necessarily affine equivariant estimator $T$ as defined above, 
$$
T^*(X)=(T_1^*(X),T_2^*(X),T_3^*(X))=
(T_1(P(X)),B(X)^{-1}(T_2(P(X)))+L(X),B(X)^{-1}T_3(P(X))\tr{(B(X)^{-1})})
$$
is affine equivariant.
\end{proposition}
\begin{proof} Let $Y=AX+a$. Observe that $B(Y)=B(X)A^{-1}$, and consequently 
$P(Y)=P(X)$. Therefore $T^*(Y)=(T^*_1(X),AT_2^*(X)+a,AT_3^*(X)\tr{A})$ as required.
\end{proof}
Note that the ERC here apply to $T^*(P(X))$, but not necessarily to $T^*(X)$. Note also that $L$ and $S$ can be any affine equivariant location and scatter estimates, for example the sample mean and covariance matrix, regardless of any knowledge of the underlying distribution or clustering.}


\addcontentsline{toc}{section}{References}
\bibliographystyle{chicago}
\bibliography{REFS-pietro,REFS-chris}%

\end{document}